\documentclass{amsart}
\usepackage[centertags]{amsmath}
\usepackage{amsfonts}
\usepackage{amssymb}
\usepackage{amscd}
\usepackage{amsthm}
\usepackage{fancyhdr}
\usepackage{graphicx}
\title{On weakly formulated Sylvester equations and applications}
\author{
Luka Grubi\v{s}i\'{c}}
\address{FernUniversit\" at in Hagen,\\
LG Mathematische Physik,\\
Feithstr. 140,\\
D-58084 Hagen}
\email{luka.grubisic@fernuni-hagen.de}
\author{
Kre\v{s}imir Veseli\'{c}}
\address{FernUniversit\" at in Hagen,\\
LG Mathematische Physik,\\
Feithstr. 140,\\
D-58084 Hagen}
\email{kresimir.veselic@fernuni-hagen.de}
\keywords{
Eigenvalues, eigenvectors,
Variational methods for eigenvalues of operators,
Perturbation theory
}

\subjclass{
65F15, 49R50, 47A55, 35Pxx}
\begin{document}
%
%  defs
%
\def\ra{{\sf R}}
\def\je{{\sf N}}
\def\I{\mathbf{I}}
\def\mA{\mathbf{A}}
\def\mM{\mathbf{M}}
\def\mU{\mathbf{U}}
\def\mX{\mathbf{X}}
\def\mV{\mathbf{V}}
\def\mW{\mathbf{W}}
\def\mY{\mathbf{Y}}
\def\mK{\mathbf{K}}
\def\mR{\mathbf{R}}
\def\mQ{\mathbf{Q}}
\def\mT{\mathbf{T}}
\def\mS{\mathbf{S}}
\def\mH{\mathbf{H}}
\def\x{\mathcal{X}}
\def\q{\mathcal{Q}}
\def\b{\mathcal{B}}
\def\g{\mathcal{G}}
\def\d{\mathcal{D}}
\def\K{\mathcal{K}}
\def\H{\mathcal{H}}
\def\vp{\mathcal{V}}
\def\wp{\mathcal{W}}
\def\lp{\mathcal{L}}
\def\EE{\mathcal{E}}
\def\SS{\mathcal{S}}
\def\R{ \mathbb{R}}
\def\C{ \mathbb{C}}
\def\N{ \mathbb{N}}

%
% environments
%
\newtheorem{theorem}{Theorem}[section]
\newtheorem{corollary}[theorem]{Corollary}
\newtheorem{lemma}[theorem]{Lemma}
\newtheorem{proposition}[theorem]{Proposition}
\theoremstyle{definition}
\newtheorem{remark}[theorem]{{\rm\textbf{Remark}}}
\newtheorem{assumption}[theorem]{{\rm \textbf{Assumption}}}
\newtheorem{example}[theorem]{{\rm \textbf{Example}}}
\newtheorem{problem}{{\rm \textbf{Problem}}}[section]
\numberwithin{equation}{section}
\renewcommand{\theequation}{\thesection.\arabic{equation}}
\newtheorem{definition}[theorem]{{\rm \textbf{Definition}}}
%
% commands
%
\newcommand{\norm}[1]{\Vert #1 \Vert}
\newcommand{\normm}[1]{\Vert #1 \Vert}
\newcommand{\ugl}[1]{\left[ #1 \right]}
\newcommand{\sk}[1]{\left( #1 \right)}
\newcommand{\dual}[1]{\left< #1 \right>}
\newcommand{\abs}[1]{\left| #1 \right|}
\newcommand{\abss}[1]{\vert #1 \vert}
\newcommand{\spa}[1]{\mbox{span}\{ #1 \}}
\newcommand{\absf}[2]{\frac{\abss{#1}}{\abss{#2}} }
\newcommand{\ine}[1]{{\mathbf{#1}}}
\newcommand{\conpo}[1]{\stackrel{#1}{\longrightarrow}}
%
% functions
%
\def\imag{{\rm i}}
\def\sinbf{\mathsf{sin}}
\def\cosbf{\mathsf{cos}}
\def\eexp{\text{e}}
\def\region{\mathcal{R}}
\def\slim{\text{{\rm s-lim}}}
\def\wlim{\text{{\rm w-lim}}}
\def\tripleb{\mid\!\mid\!\mid}

%
% Faks (kuci)
%

%\def\Weidmann_I{WeidmannBuch}
\def\Weidmann_I{Weidmann_I}

\begin{abstract}
We use a ``weakly formulated'' Sylvester equation $$A^{1/2}TM^{-1/2}-A^{-1/2}TM^{1/2}=F$$
to obtain new bounds for the rotation of spectral
subspaces of a nonnegative selfadjoint operator in a Hilbert space.
Our bound extends the known results of Davis and Kahan.
Another application is a bound for the square root of a positive selfadjoint
operator which extends the known rule: ``The relative
error in the square root is bounded by the one half of the relative error in
the radicand''. Both bounds are illustrated on differential operators which are
defined via quadratic forms.
\end{abstract}

%\begin{keywords}
%Eigenvalues,
%Estimation of eigenvalues, upper and lower bounds,
%Variational methods for eigenvalues of operators
%\end{keywords}
%
%\begin{AMS}
%65F15, 34L15, 35P15, 49R50
%\end{AMS}
\maketitle
\pagestyle{myheadings}
\thispagestyle{plain}
\markboth{L. Grubi\v{s}i\'{c} and Kre\v{s}imir Veseli\'{c} }{L. Grubi\v{s}i\'{c} and
Kre\v{s}imir Veseli\'{c}}

\section{Preliminaries}
In this work we will study properties of nonnegative selfadjoint operators in a Hilbert space
which are close in the sense
of the inequality
\begin{equation}\label{eq:prva_ocjena}
|h(\phi,\psi) - m(\phi,\psi)| \leq \eta \sqrt{h[\phi]m[\psi]}
\end{equation}
where the sesquilinear forms \(h, m\) belong to the operators
\(\mH,\mM\) respectively and $m[\psi]=m(\psi,\psi)$, $h[\phi]=h(\phi,\phi)$.

In the first part of the paper we show that (\ref{eq:prva_ocjena})  implies an
estimate of the separation between ``matching''
eigensubspaces of $\mH$ and $\mA$. To be more precise
one of the typical situations is: Let
\begin{align}\label{eigenvaluesH}
0&\leq\lambda_1(\mH)\leq\lambda_2(\mH)\leq
\cdots\leq\lambda_n(\mH)<D<\lambda_{n+1}(\mH)\leq\cdots\\
0&\leq\lambda_1(\mM)\leq\lambda_2(\mM)\leq
\cdots\leq\lambda_n(\mM)<D<\lambda_{n+1}(\mM)\leq\cdots
\label{eigenvaluesM}
\end{align}
be the eigenvalues of the operators $\mH$ and $\mM$ which satisfy (\ref{eq:prva_ocjena})
then
$$
\|E_{\mH}(D)-E_{\mM}(D)\|\leq\min\Big\{\frac{\sqrt{D\lambda_n(\mH)}}{D-\lambda_n(\mH)},
\frac{\sqrt{D\lambda_n(\mM)}}{D-\lambda_n(\mM)}\Big\}~\eta.
$$
Such an estimate\footnote{For recent estimates of the separation
between eigensubspaces see \cite{Kostrykin05}.} was implicit in \cite{GruVes02}.  We then generalize
this inequality to hold both for the operator norm $\|\cdot\|$ and the Hilbert--Schmidt norm
 $\tripleb\cdot \tripleb_{HS}$. We also allow that
%\footnote{When we say that some
%orthogonal projection $P$ is infinite dimensional we mean that
%$\ra(P)$ is infinite dimensional subspace of the environment Hilbert space.}
$E_\mH(D)$ and $E_\mM(D)$ be possibly infinite dimensional\footnote{We assume that
$E(\cdot)$ is right continuous.}.

In the second part of the paper we establish estimates for
a perturbation of the square root of a positive operator.
It will be shown that the inequality (\ref{eq:prva_ocjena})
implies
\[
|h_2(\phi,\psi) - m_2(\phi,\psi)|
\leq \frac{\eta}{2} \sqrt{h_2[\phi]m_2[\psi]},
\]
where the sesquilinear forms \(h_2, m_2\) belong to the operators
\(\mH^{1/2}, \mM^{1/2}\), respectively. This will show that it is meaningful to
consider weakly formulated Sylvester equations where all the coefficient operators
are unbounded, cf. \ref{weakS}.

Both of this problems will be solved through a study
of the weak Sylvester equation, which reads formally
\begin{equation}\label{weakS}
\mH X- X \mM=\mH^{1/2}F\mM^{1/2}.
\end{equation}
These two case studies represent two different classes of additional assumptions which
have to be imposed on the coefficient operators $\mH$, $\mM$ and $F$ in order that
(\ref{weakS}) defines a meaningful operator $X$.

The main novelty (and contribution) of this work is that we present an abstract study of the operator
equation (\ref{weakS}) in the case when only $F$ is a \textit{bona fide} operator.
The expression $\mH^{1/2}F\mM^{1/2}$ need not possess an operator representation.
In comparison, $\mH^{1/2}F\mM^{1/2}$ was always a bounded operator for
the Sylvester equations which were studied in
\cite{BathDavMcInt85,DavisKahan70,Ren-CangStructured}.

Our first main result, contained in Theorem \ref{prvo2:t_slabisylv} below, extends
our previous result from \cite{GruVes02} in various ways. In particular, we allow
the perturbed projection to be infinite dimensional. In the proof we also
overcome a technical error contained in \cite{GruVes02}. We then extend
this result to the case of other unitary invariant operator norms\footnote{
Also called ``cross-norms'' in the terminology of \cite{Kato76} or ``symmetric norms'' in the terminology
of \cite{GohbergKrein,SimonTrace}.} Particular
attention is paid to the Hilbert--Schmidt norm because of its
possible importance in applications. This special case is handled by
another technique which allows an arbitrary interlacing of
the involved spectra.

The main object in this work shall be a closed nonnegative symmetric
form in a Hilbert space.
When dealing with symmetric forms in a Hilbert space, we
shall follow the terminology of Kato, cf. \cite{Kato76}.
For reader's convenience we now give
definitions of some terms that will frequently be used, cf. \cite{Faris75,Kato76}.

\begin{definition}
Let $h$ be a positive definite form in $\H$. A sesquilinear form $a$, which need not
be closed, is said to be \textit{$h$-bounded},\index{form!$h$-bounded form} if $\q(h)\subset\q(a)$ and there exists $\eta\geq 0$
$$
|a[u]|\leq \eta h[u] \qquad u\in\q(h).
$$
%The greatest lower bound of all $b_h(h')$ is called the $h$-bound of $h'$.
\end{definition}
If $h$ is positive definite the space $(\q(h),h)$ can be considered as
a Hilbert space. The form $a$, which is $h$-bounded,
defines a bounded operator on the space $(\q(h),h)$.

\begin{definition}
A bounded operator $A:\H\to\mathcal{U}$ is called \textit{degenerate}\index{operator!degenerate}
if $\ra(A)$ is finite dimensional.
\end{definition}
\begin{definition}\label{prvo:d_komutira}
If $\mH$ is a self adjoint operator and $P$ a projection, to say that
$P$ commutes with $\mH$ means that $u\in\d(\mH)$ implies $Pu\in\d(\mH)$
and
$$
\mH Pu=P\mH u,\qquad u\in\d(\mH).
$$
\end{definition}

\begin{definition}
Let $\mH$ and $\mM$ be nonnegative operators. We define the \textit{order relation}
\index{operator!the order relation}
$\leq$ between the nonnegative operators by saying that
$$
\mM\leq\mH
$$
if and only if $\q(\mH)\subset\q(\mM)$ and
$$
\|\mM^{1/2} u\|\leq\|\mH^{1/2} u\|,\qquad u\in\q(\mH),
$$
or equivalently
$$m[u]\leq h[u], \qquad u\in\q(h),$$
when $m$ and $h$ are nonnegative forms defined by the operators $\mM$ and
$\mH$ and $\mM\leq\mH$.
\end{definition}

A main principle we shall use to develop the perturbation theory will be
the \textit{monotonicity of the spectrum}\index{spectrum!monotonicity of the spectrum}
with regard to the order relation
between nonnegative operators. This principle can be expressed in many ways.
The relevant results, which are scattered over the monographs \cite{Faris75,Kato76}, are summed
up in the following theorem, see also \cite[Corollary A.1]{Levendorskii}.

\begin{theorem}\label{prvo:t_monotonost}
Let $\mM=\int\lambda~\text{d}E_\mM(\lambda)$ and $\mH=\int\lambda~\text{d}E_\mH(\lambda)$
be nonnegative operators in $\H$ and let $\mM\leq\mH$. Let the eigenvalues
of $\mH$ and $\mM$ be as in (\ref{eigenvaluesH}) and (\ref{eigenvaluesM}) then
\begin{enumerate}
\item $\lambda_e(\mM)\leq\lambda_e(\mH)$
\item ${\sf dim} ~E_\mH(\gamma)\leq {\sf dim} ~E_\mM(\gamma)$, for every $\gamma\in\R$
\item $\lambda_k(\mM)\leq\lambda_k(\mH),\qquad k=1,2,\cdots$ .
\end{enumerate}
The infimum of the essential spectrum of some operator $\mH$ is denoted by $\lambda_e(\mH)$.
\end{theorem}

With this theorem in hand we review spectral properties of operators $\mH$ and $\mM$,
for which there exists $0\leq\varepsilon<1$ such that
\begin{equation}\label{prvo:e_prvizakjN}
(1-\varepsilon)m[u]\leq h[u]\leq(1+\varepsilon)m[u],\qquad u\in\q:=\q(h)=\q(m).
\end{equation}
Let us assume $h[u]>0$ then $m[u]>0$ and
\begin{equation}\label{eq:brojevi}
(1-\frac{\varepsilon}{1-\varepsilon})h[u]\leq m[u]\leq (1+\frac{\varepsilon}{1-\varepsilon})h[u].
\end{equation}
Inequality (\ref{prvo:e_prvizakjN}) implies that $\je(\mH)=\je(\mM)$, so (\ref{eq:brojevi})
holds for all $u\in\q$.
%\begin{equation}\label{prvo:e_drugizakjN}
%(1-\frac{\varepsilon}{1-\varepsilon})h[u]\leq m[u]\leq (1+\frac{\varepsilon}{1-\varepsilon})h[u],
%\qquad u\in\q.
%\end{equation}

\begin{lemma}\label{lemma:comb_1}
Let $m$ and $h$ be nonnegative forms such that $\lambda_e(\mM)>0$ and $\lambda_e(\mH)>0$ and
let (\ref{prvo:e_prvizakjN}) hold. Then
\begin{align}\label{eq.rel_1}
|\lambda_i(\mH)-\lambda_i(\mM)|&\leq\varepsilon\lambda_i(\mM)\\
|\lambda_i(\mH)-\lambda_i(\mM)|&\leq\frac{\varepsilon}{1-\varepsilon}\lambda_i(\mH)
\label{eq.rel_2}
\end{align}
$\lambda_i(\mH)$ and $\lambda_i(\mM)$ are as in (\ref{eigenvaluesH}) and (\ref{eigenvaluesM}).
Assume that $\lambda_{i-1}(\mH)<\lambda_i(\mH)<\lambda_{i+1}(\mH)$ and
\begin{equation}\label{eq.uvijet12}
\frac{\varepsilon}{1-\varepsilon}<\max\Big\{
\frac{\lambda_{i+1}(\mH)-\lambda_i(\mH)}{\lambda_{i+1}(\mH)+\lambda_i(\mH)},
\frac{\lambda_{i}(\mH)-\lambda_{i-1}(\mH)}{\lambda_{i}(\mH)+\lambda_{i-1}(\mH)},1\Big\}
\end{equation}
then
\begin{equation}\label{eq.prva_dok}
\min_{\lambda_j(\mM)}\frac{|\lambda_i(\mH)-\lambda_j(\mM)|}{\lambda_i(\mH)}=
\frac{|\lambda_i(\mH)-\lambda_i(\mM)|}{\lambda_i(\mH)}<1.
\end{equation}
If $\lambda_{i-1}(\mH)<\lambda_i(\mH)=\cdots=\lambda_{i+n-1}(\mH)<\lambda_{i+n}(\mH)$
and
\begin{equation}\label{eq.uvijet123}
\frac{\varepsilon}{1-\varepsilon}<\max\Big\{
\frac{\lambda_{i+n}(\mH)-\lambda_i(\mH)}{\lambda_{i+n}(\mH)+\lambda_i(\mH)},
\frac{\lambda_{i}(\mH)-\lambda_{i-1}(\mH)}{\lambda_{i}(\mH)+\lambda_{i-1}(\mH)},1\Big\}
\end{equation}
then
\begin{align}\label{eq.druga_dok_1}
\underset{j\in\N}{\text{{\rm argmin}}} \frac{|\lambda_{i-1}(\mH)-\lambda_j(\mM)|}{\lambda_{i-1}(\mH)}&\leq i-1\\
\underset{j\in\N}{\text{{\rm argmin}}} \frac{|\lambda_{i+n}(\mH)-\lambda_j(\mM)|}{\lambda_{i+n}(\mH)}&\geq i+n.
\label{eq.druga_dok_2}
\end{align}
\end{lemma}
\begin{proof}
Estimates (\ref{eq.rel_1})--(\ref{eq.rel_2}) are a
 consequence of (\ref{prvo:e_prvizakjN})--(\ref{eq:brojevi}) and Theorem
\ref{prvo:t_monotonost}. The rest of the theorem follows from a proof which
analogous to the proof of \cite[Theorem 4.16]{Gru03_3}. We repeat the argument in
this new setting.

Let $i\ne j$ then
\begin{align*}
\frac{|\lambda_i(\mH)-\lambda_j(\mM)|}{\lambda_i(\mH)}&\geq
\frac{|\lambda_i(\mH)-\lambda_j(\mH)|}{\lambda_i(\mH)+\lambda_j(\mH)}\frac{\lambda_i(\mH)+\lambda_j(\mH)}{\lambda_i(\mH)}
-\frac{|\lambda_j(\mH)-\lambda_j(\mM)|}{
\lambda_j(\mH)}\frac{\lambda_j(\mH)}{\lambda_i(\mH)}\\
&\geq\gamma \Big(1+\frac{\lambda_j(\mH)}{\lambda_i(\mH)}\Big)-
\frac{\varepsilon}{1-\varepsilon} \frac{\lambda_j(\mH)}{\lambda_i(\mH)}>\gamma\\
& >\frac{|\lambda_i(\mH)-\lambda_i(\mM)|}{\lambda_i(\mH)}.
\end{align*}
With this we have established (\ref{eq.prva_dok}). (\ref{eq.druga_dok_1})--(\ref{eq.druga_dok_2})
are a way to state (\ref{eq.prva_dok}) in a presence of a multiple eigenvalue $\lambda_i(\mH)$.
The proof follows by a repetition of the previous argument for $j\geq i$ and $j\leq i+n-1$.
For instance, we establish (\ref{eq.druga_dok_1}) by proving
$$
\frac{|\lambda_{i-1}(\mH)-\lambda_j(\mM)|}{\lambda_{i-1}(\mH)}>
\frac{|\lambda_{i-1}(\mH)-\lambda_{i-1}(\mM)|}{\lambda_{i-1}(\mH)}
$$
for all $j\geq i$.
\end{proof}

\begin{remark}
The significance of this Lemma is that it detects which spectral subspaces should be compared.
When we were comparing discrete eigenvalues, the order relation between the real numbers (eigenvalues)
solved this problem automatically. For spectral subspaces we need to assume more than
(\ref{prvo:e_prvizakjN}) in order to be able to construct meaningful estimates. Assumptions
(\ref{eq.uvijet12}) and (\ref{eq.uvijet123}) show how much more we (will) assume.
\end{remark}

Next we show that (\ref{prvo:e_prvizakjN}) implies (\ref{eq:prva_ocjena}) with
$\eta=\varepsilon(1-\varepsilon)^{-1/2}$.
To establish this claim we need a notion of a
pseudo inverse of a closed operator. A definition from \cite{WeidmannScand84} will be used. The
\textit{pseudo inverse}\index{generalized inverse!pseudo inverse}
of a self adjoint operator $\mH$ is the self adjoint operator $\mH^\dagger$ defined by
\begin{align*}
\d(\mH^\dagger)&=\ra(\mH)\oplus\d(\mH)^\perp,\\
\mH^\dagger(u+v)&=\mH^{-1}u,\qquad u\in\ra(\mH),~g\in\d(\mH)^\perp.
\end{align*}
It follows that $\mH^\dagger=\mH^{-1}$ in $\overline{\ra(\mH)}$. Note that we did not assume
$\mH^\dagger$ to be bounded or densely defined. The operator $\mH^\dagger$ will be bounded if an
only if $\ra(\mH)$ is closed in $\H$, see \cite{Nashed76}. The operator $\mH^\dagger$ could
have also been defined by the spectral calculus, since
$$
\mH^\dagger=f(\mH),\qquad f(\lambda)=\begin{cases}0,& \lambda=0,\\\frac{1}{\lambda},&\lambda\ne 0.
\end{cases}
$$
In \cite{WeidmannScand84}
Weidmann has given a short survey of the properties
of the pseudo inverse of a nondensely defined operator $\mH$.
In particular, let $\mH_1$ and $\mH_2$ be
two nonnegative operators in $\overline{\d(\mH_1)}$ and $\overline{\d(\mH_2)}$ respectively then
\begin{equation}\label{e_inveruredjaj}
\|\mH^{1/2}_1u\|\leq\|\mH^{1/2}_2u\|\Leftrightarrow\|\mH^{1/2\dagger}_2u\|\leq\|\mH^{1/2\dagger}_1u\|.
\end{equation}
Analogously, let $h_1$ and $h_2$ be two closed, not necessarily densely defined, positive definite forms and let
$\mH_1$ and $\mH_2$ be the self adjoint operators defined by $h_1$ and $h_2$ in $\overline{\q(h_1)}$
and $\overline{\q(h_2)}$. We say $h_1\leq h_2$ when $\q(h_2)\subset\q(h_1)$ and
\begin{equation}\label{e_uredjaj}
h_1[u]=\|\mH_1^{1/2}u\|^2\leq h_2[u]=\|\mH_2^{1/2}u\|^2,\qquad u\in\q(h_2).
\end{equation}
Equivalently, we write $\mH_1\leq\mH_2$ when $h_1\leq h_2$. Now, we can write the fact (\ref{e_inveruredjaj})
as
\begin{equation}\label{e_radniuredjaj}
\mH_1\leq\mH_2\Longleftrightarrow\mH_2^{\dagger}\leq\mH^{\dagger}_1.
\end{equation}

In one
point we will depart from the conventions in \cite{Kato76}.
\begin{definition}
A nonnegative form $$h(u,v)=(\mH^{1/2}u, \mH^{1/2}v)$$
will be called \emph{nonnegative definite}\index{form!nonnegative definite}
when  $\mH^\dagger$ is bounded. Analogously, the nonnegative operator $\mH$ such that
$\mH^\dagger$ is bounded will also be called
\emph{nonnegative definite}.

\end{definition}
In the sequel we establish a connection between (\ref{prvo:e_prvizakjN}) and
(\ref{eq:prva_ocjena}) when $h$ and $m$ are nonnegative definite forms.

\begin{lemma}\label{prvo1:l_mathves}
Let \(\mH\) and \(\mM\) be nonnegative definite operators in a Hilbert
space \(\H\) such that (\ref{prvo:e_prvizakjN}) holds
for $0\leq\varepsilon <1$.
Let
\begin{equation}\label{Tgeneral1}
S=\mH^{1/2}\mM^{\dagger 1/2}-\overline{\mH^{\dagger 1/2}\mM^{1/2}}
\end{equation}
then $S$ is bounded and
\begin{equation}\label{prvo1:e_mathves}
|(\psi,S \phi)|\leq\frac{\varepsilon}{\sqrt{1-\varepsilon}}\|\psi\|\|\phi\|.
\end{equation}

\end{lemma}
\begin{proof}
The closed graph theorem implies that
the operator
$$
S=\mH^{1/2}\mM^{\dagger 1/2}-\overline{\mH^{\dagger 1/2}\mM^{1/2}}
$$
is bounded. Also, $\je(\mH)=\je(\mM)=\je(S)$ and $P_{\je(S)}$ commutes with $S$.
It is sufficient to prove the estimate for $x, y\in\ra(\mH)$.
The assumption (\ref{prvo:e_prvizakjN})
gives
$$
|\big(h-m\big)(\mH^{\dagger 1/2} x,\mM^{\dagger 1/2} y)|\leq\varepsilon\norm{ y}~
m[\mH^{\dagger 1/2} x]^{1/2}.
$$
Analogously, (\ref{prvo:e_prvizakjN}) implies
\begin{equation}\label{prvo:e_ocjenaZaUpotr}
\|\mM^{1/2}\mH^{\dagger 1/2}\|\leq\frac{1}{\sqrt{1-\varepsilon}}~.
\end{equation}
Altogether, the estimate (\ref{prvo1:e_mathves}) follows.
\end{proof}
Now, we rewrite the conclusion of this lemma in the symmetric form setting.
The result is given in the form of a proposition which we present without proof.
\begin{proposition}\label{prvo1:t_dekompozicija}
Let \(m\) and \(h\) be nonnegative definite forms and let there exist  $0\leq\varepsilon<1$ such
that (\ref{prvo:e_prvizakjN}) holds then $\je(\mH)=\je(\mH)$ and
$$
|h(u,v) - m(u,v)| \leq \frac{\varepsilon}{\sqrt{1-\varepsilon}}\sqrt{h[u]m[v]}.
$$
\end{proposition}
When we only know that $h$ and $m$ satisfy (\ref{eq:prva_ocjena}) then we can establish
a similar result about $\je(\mH)$ and $\je(\mM)$.
\begin{proposition}\label{lema:sim}
Let $m$ and $h$ be nonnegative definite forms such that (\ref{eq:prva_ocjena}) holds
then \begin{align*}
S&=\mH^{1/2}\mM^{\dagger 1/2}-\overline{\mH^{\dagger 1/2}\mM^{1/2}}\\
S^*&=\overline{\mM^{\dagger 1/2}\mH^{1/2}}-\mM^{1/2}\mH^{\dagger 1/2}
\end{align*}
are bounded operators and $\|S^*\|=\|S\|\leq\eta$. Furthermore,
$\je(\mH)=\je(\mM)$ and a fortiori $\ra(\mH)=\ra(\mM)$.
\end{proposition}

 The operator $S$ has a special structure.
Assume $\mM u=\mu u$ and $\mH v=\lambda v$,
then
\begin{align}
\nonumber (v, Su)&=\lambda^{1/2}(v, u)\mu^{1/2} -\lambda^{-1/2}(v, u)\mu^{1/2}\\
&=\frac{\lambda-\mu}{\sqrt{\lambda\mu}}(v,u)\label{prvo:eq_strukruraSa}~.
\end{align}
The equation (\ref{prvo:eq_strukruraSa}) suggests the distance function
$$
\frac{|\lambda-\mu|}{\sqrt{\lambda\mu}}
$$
which measures the distance between the eigenvalues of operators $\mH$ and $\mM$.
We state this result as the following corollary..
\begin{corollary}
Let $\mM u=\mu u$, $\|u\|=1$ and $\mH v=\lambda v$, $\|v\|=1$
and let $S$ be as in Proposition \ref{lema:sim} then
$$
\frac{|\lambda-\mu|}{\sqrt{\lambda\mu}}\leq\frac{\eta}{|(u,v)|}.
$$
\end{corollary}

Our theory is designed to be directly applicable to differential
operators given in a weak form. This will enable us to obtain estimates
for the difference between the spectral projections of the operators to
which the theory of \cite{BathDavMcInt85,DavisKahan70} does not apply, see Example
\ref{example:measur} below.

\section{Weak Sylvester equation}\label{Sylvester equation}

Let us outline the general picture.
We have an unbounded positive definite operator
$\mathbf{A}$ and a bounded positive definite operator $M$. They are defined
in, possibly, different subspaces of the environment Hilbert space $\H$.
Thus, $\H_M=\ra(M)$ is (of necessity) a closed subspace of $\H$ and
likewise $$\overline{\d(\mA^{1/2})}^{~_\H}=\ra(\mA^{1/2})=\H_\mA.$$
Let the
bounded operator $F:\H_M\to\H_\mA$ be given, then
we are looking for the bounded operator $T:\H_M\to\H_\mA$ such that
\begin{equation}\label{drugo:weak_sylvester}
(\mA^{1/2}  v,TM^{-1/2} u)-(\mA^{- 1/2}v,TM^{1/2} u)=(v, Fu)~,
\qquad v\in\d(\mA^{1/2}),~ u\in \H_M.
\end{equation}
 Formally, we say that $T$ solves the equation
\begin{equation}\label{prvo2:e_formalnaS}
\mA T- TM=\mA^{1/2}FM^{1/2}.
\end{equation}
Here $G=\mA^{1/2}FM^{1/2}$ is naturally only a formal expression and does
not represent a \textit{bona fide} operator. In the case in which $G$
be a \textit{bona fide} operator equation (\ref{prvo2:e_formalnaS}) becomes the
rigorous equation
$$
\mA T- TM=G,
$$
called the (standard) \textit{Sylvester equation}, cf. \cite{BathDavMcInt85,DavisKahan70}.
The case when $\mA$ and $M$ are finite matrices
has been considered in \cite{Ren-CangStructured} where (\ref{prvo2:e_formalnaS})
was called the structured Sylvester equation.

We call the relation (\ref{drugo:weak_sylvester}) the \textit{weak Sylvester
equation}.
It represents a generalization of the concept of
the structured Sylvester equation (\ref{prvo2:e_formalnaS}) from
finite matrix setting to unbounded operator setting.
The following theorem slightly generalizes the corresponding
result from the joint paper \cite{GruVes02} and corrects a technical
glitch in one of the proofs.

\begin{theorem}\label{prvo2:t_slabisylv}
Let $\mathbf{A}$ and $M$ be positive definite operators in $\H_{\mathbf{A}}$
and $\H_M$, respectively and let $F$ be a bounded operator from
$\H_M$ into $\ra(\mathbf{A}^{1/2})=\H_{\mathbf{A}}$. If $M$ is bounded and
\begin{equation}\label{prvo2:e_dihotomija2}
\norm{M}<\frac{1}{\normm{\mathbf{A}^{-1}}}
\end{equation}
then the weakly formulated Sylvester equation
\begin{equation}\label{prvo2:e_slabaS1}
\sk{\mathbf{A}^{1/2} v, TM^{-1/2}u}-\sk{v,\mathbf{A}^{-1/2}TM^{1/2}u}=\sk{v, F u}
\end{equation}
has a unique solution $T$, given by $\tau(v,u)=(v, Tu)$ and
\begin{equation}\label{prvo2:e_integralnaR}
\tau(v, u)=-\frac{1}{2\pi} \int_{-\infty}^{\infty}
(\mathbf{A}^{1/2}v, (\mathbf{A}-{\rm i}\zeta-d)^{-1}F(M-{\rm i}\zeta -d)^{-1}M^{1/2}u)d\zeta,
\end{equation}
where $d$ is any number satisfying
\begin{equation}\label{prvo2:e_dihotomija}
\norm{M}<d<\frac{1}{\normm{\mathbf{A}^{-1}}}\;.
\end{equation}
%provided that
%\begin{equation}
%\norm{M}<\frac{1}{\normm{\mathbf{A}^{-1}}}.
%\end{equation}
\end{theorem}
\proof
The uniqueness means that
\begin{equation}\label{prvo2:e_dvasest}
\sk{\mathbf{A}^{1/2} v,WM^{-1/2}u}-\sk{v,\mathbf{A}^{-1/2}WM^{1/2}u}=0,
\end{equation}
for $u\in\H_M$, $v\in \d(\mathbf{A}^{1/2})$, has the only bounded solution $W=0$.
Let $$E_n=\int_0^n d~E_{\mA^{1/2}}(\lambda),$$ then in particular
$$
\sk{\mathbf{A}^{1/2} v,E_nWM^{-1/2}u}-\sk{v,\mathbf{A}^{-1/2}E_nWM^{1/2}u}=0,
$$
for $u\in\H_M$, $v\in \d(\mathbf{A}^{1/2})\cap E_n\H$.
Define the cut--off function
$$
f_n(x)=\left\{\begin{array}{ll}
x,& D\leq x\leq n\\
n,& n \leq x\end{array}\right.
$$
with $D=1/\|\mathbf{A}^{-1}\|$. The operator $f_n(\mathbf{A}^{1/2})$ is bicontinuous and
\begin{equation}\label{drugo:standardsylv}
f_n(\mathbf{A}^{1/2})E_n WM^{-1/2}-f_n(\mathbf{A}^{1/2})^{-1}E_n WM^{1/2}=0.
\end{equation}
Since $f_n(\mathbf{A}^{1/2})$ and $M^{1/2}$ are bounded and positive definite operators,
the standard Sylvester equation
(\ref{drugo:standardsylv})
has the unique solution
\begin{equation}\label{prvo:e_E_nW}
E_n W=0,\qquad n\in\N~.
\end{equation}
This is a consequence of the standard theory of the Sylvester equation with
bounded coefficients, see \cite{BathDavMcInt85,DavisKahan70}.
The statement (\ref{prvo:e_E_nW}) implies $W=0$.

Now for the existence.  We use
the spectral integral $\mA=\int\lambda \;d E(\lambda)$ to compute
\begin{align}
\nonumber \int^\infty_{-\infty}\|(\mA+i\zeta-d)^{-1}\mA^{1/2}v\|^2~&d\zeta=
 \int^\infty_{-\infty} (\mA^{1/2}v,  \big|\mA-i\zeta-d\big|^{-2}\mA v)~d\zeta\\
\nonumber
&=\int^\infty_{-\infty}~d\zeta\int_D^\infty \frac{\lambda~d(E(\lambda)\mA^{1/2}v,\mA^{1/2}v)}
{(\lambda-d)^2+\zeta^2}\\
\nonumber &=\int_D^\infty \lambda~d(E(\lambda)\mA^{1/2}v, \mA^{1/2}v)\int^\infty_{-\infty}\frac{d\zeta}{
(\lambda-d)^2+\zeta^2}\\
\nonumber &=\int^\infty_D\frac{\pi \lambda ~d(E(\lambda)\mA^{1/2}v, \mA^{1/2}v)}{\lambda-d}\\
&=\pi(\mA (\mA-d)^{-1} v, v).\label{prvo:integral1}
\end{align}
Analogously, one establishes
\begin{equation}\label{prvo:integral2}
\int^\infty_{-\infty}\|(M-i\zeta-d)^{-1}M^{1/2}u\|^2~d\zeta= \pi
(M (d-M)^{-1} u, u).
\end{equation}
The convergence of these integrals justifies the following computation. Set
$$
\tau(v,u)=-\frac{1}{2\pi} \int_{-\infty}^{\infty}
(\mathbf{A}^{1/2}v, (\mathbf{A}-{\rm i}\zeta-d)^{-1}F(M-{\rm i}\zeta -d)^{-1}M^{1/2}u)d\zeta
$$
and then compute using (\ref{prvo:integral1}) and (\ref{prvo:integral2})
\begin{align}
\nonumber|\tau(v, u)|^2&=\frac{1}{(2\pi)^2}\Big[\int^\infty_{-\infty}
((\mathbf{A}+{\rm i}\zeta-d)^{-1}\mathbf{A}^{1/2}v, F(M-{\rm i}\zeta -d)^{-1}M^{1/2}u)d\zeta\Big]^2\\
\nonumber&\leq\frac{\|F\|^2}{(2\pi)^2}
\Big[\int^\infty_{-\infty}\|(\mA+i\zeta-d)^{-1}\mA^{1/2}v\|\;
\|(M-i\zeta-d)^{-1}M^{1/2}u\| d\zeta\Big]^2\\
&\leq\frac{\|F\|^2}{4}(\mA (\mA-d)^{-1} v, v)(M (d-M)^{-1} u, u).\label{eq:prerequisite}
\end{align}
This in turn implies that the operator
$$
\tau(v, u)=(v, Tu)
$$
is a bounded operator and also gives the meaning to the formula (\ref{prvo2:e_integralnaR}).

Now we will prove that this $T$ satisfies the equation (\ref{prvo2:e_slabaS1}). Note that
$$
\mA(\mA-\rho-d)^{-1}=\I+(\rho+d)(\mA-\rho-d)^{-1},\qquad \rho\not\in\sigma(\mA)
$$
and then take $v\in\d(\mA)$ to compute
\begin{align*}
(\mA^{1/2}v, TM^{-1/2}u) &- (\mA^{-1/2}v, TM^{1/2}u)=\\
&=-\frac{1}{2\pi}\Big[\int^\infty_{-\infty}
(\mathbf{A}v, (\mathbf{A}-{\rm i}\zeta-d)^{-1}F(M-{\rm i}\zeta -d)^{-1}u)~d\zeta\\
&\phantom{\Big[\int^\infty_{-\infty}}-\int^\infty_{-\infty}
(v, (\mathbf{A}-{\rm i}\zeta-d)^{-1}F(M-{\rm i}\zeta -d)^{-1}M u)~d\zeta\Big]\\
&=-\frac{1}{2\pi}\Big[
{\rm \text{v.p.}}\int^\infty_{-\infty}(v, F(M-i\zeta-d)^{-1} u)~d\zeta\\
&\phantom{\Big[\int^\infty_{-\infty}}+\int^\infty_{-\infty}
(i\zeta + d)((\mA-i\zeta-d)^{-1} v, F(M-i\zeta-d)^{-1} u)~d\zeta\\
&\phantom{\Big[\int^\infty_{-\infty}}-\int^\infty_{-\infty}
(i\zeta + d)((\mA-i\zeta-d)^{-1} v, F(M-i\zeta-d)^{-1} u)~d\zeta\\
&\phantom{\Big[\int^\infty_{-\infty}}-{\rm \text{v.p.}}\int^\infty_{-\infty}
( (\mA-i\zeta-d)^{-1} v, F u)~d\zeta
\Big]\\
&=(v, Fu).
\end{align*}
By a usual density argument we conclude that
the operator $T$ satisfies (\ref{prvo2:e_slabaS1}).\qed
\begin{theorem}\label{cor:unitary_invariant}
Let $\mA$, $M$ and $F$ be as in Theorem \ref{prvo2:t_slabisylv} then
$$
\| T \| \leq
\sqrt{\frac{D\normm{M}}{(D-d)(d-\normm{M})}}~ \frac{\| F \|}{2}
$$
for any $\|M\|<d<D$. The optimal $d$ is $d=(\|M\|+D)/2$ and then we obtain
\begin{equation}\label{eq:invers_Sylvester}
\| T \| \leq
\frac{\sqrt{D\normm{M}}}{(D-\normm{M})}~ \| F \|
\end{equation}
\end{theorem}
\begin{proof}
Estimate (\ref{eq:prerequisite}) yields
$$
\| T \| \leq \frac{\|F\|}{2}\|\mA(\mA-d)^{-1}\|~\|M(d-M)^{-1}\|\leq
\frac{\|F\|}{2}\sqrt{\frac{D\normm{M}}{(D-d)(d-\normm{M})}}
$$
This in turn implies the desired estimate. The optimality of the
$d=(\|M\|+D)/2$ can now be checked by a direct computation.
\end{proof}

\begin{remark}
In fact, we will se
that the estimate of Theorem \ref{cor:unitary_invariant} is optimal in the following sense.
Let us consider the equation (\ref{prvo2:e_slabaS1}) in another light.
Theorem \ref{prvo2:t_slabisylv} gives a set of conditions when the equation (\ref{prvo2:t_slabisylv})
has a unique solution. Theorem \ref{cor:unitary_invariant} then provides
us with an estimate of this solution.

Since for given $F$, under the conditions of Theorem \ref{prvo2:t_slabisylv}, there exists
the unique $T$ such that (\ref{prvo2:e_slabaS1}) holds, we can define the so called
``Sylvester operator'' which associates the solution $T$ to every operator $F$. The estimate
(\ref{eq:invers_Sylvester}) is then an estimate of the norm of the inverse of such an operator.

The bound (\ref{eq:invers_Sylvester}) is sharp in this sense as shows the following
example. Let $M$ and $\mA$ be such that
$$
Mq=\|M\|q,\qquad \mA p=Dp,
$$
for $p$ and $q$ one dimensional projections and let $F=pq$.

Then (\ref{prvo2:e_slabaS1}) is obviously satisfied by
$$
T=\frac{\sqrt{D\|M\|}}{D-\|M\|}pq.
$$
\end{remark}

\subsection{Allowing for a more general relation between $\sigma(M)$
and $\sigma(\mA)$.}
An analogue of Theorem \ref{prvo2:t_slabisylv} holds, if the assumption (\ref{prvo2:e_dihotomija})
is replaced by a more general one, namely that the interval
$$\left[\normm{M^{-1}}^{-1},\normm{M}\right]$$ be contained in the resolvent set
of the operator $\mA$.  We omit
the proof of the following result.
%\begin{figure}
%\begin{center}
%\input SAVE.pstex_t
%\caption{The spectral gaps}\label{prvo2:f_rupa}
%\end{center}
%\end{figure}
\begin{theorem}\label{prvo:c_Sylvester}
Let the operators $\mA$, $M$ and $F$ be as in Theorem \ref{prvo2:t_slabisylv},
and let their spectra be arranged so that
$$\sigma(\mA)\subset\big<0,D_{-}\big]\cup\big[D_{+}, \infty\big>,$$
where $0<D_{-}<\normm{M^{-1}}^{-1}$, $\normm{M}<D_{+}$.
Then (in the sense of (\ref{prvo2:e_integralnaR}))
\begin{eqnarray*}
T&=&-\frac{1}{2\pi} \int_{-\infty}^{\infty}
\mA^{1/2}(\mA-{\rm i}\zeta-d)^{-1}F (M-{\rm i}\zeta -d)^{-1}M^{1/2}d\zeta\\
&&\quad+
\frac{1}{2\pi} \int_{-\infty}^{\infty}
\mA^{1/2}(\mA-{\rm i}\zeta-g)^{-1}F (M-{\rm i}\zeta -g)^{-1}M^{1/2}d\zeta,
\end{eqnarray*}
where $D_{-}<g<\|M^{-1}\|^{-1}$ and $\|M\|<d<D_{+}$,
is the solution of the weak Sylvester equation (\ref{prvo2:e_slabaS1}). We also have the estimate
$$
\|T\|\leq\Big(\frac{\sqrt{\|M^{-1}\|^{-1}D_{-}}}{\|M^{-1}\|^{-1}-D_{-}}+
\frac{\sqrt{D_{+}\|M\|}}{D_{+}-\|M\|}\Big)\|F\|.
$$
\end{theorem}

\subsection{Estimates in the Hilbert--Schmidt norm}

%In our considerations of the Sylvester equation we have seen that
%the operator $Q$, appearing on the right hand side of the equation,
%determines the type of the solution.
%The properties of the operator $Q$ and the solution
%$T$ are connected in a much deeper sense than is
%signaled by Corollary \ref{prvo:degenerate_ideal}.

%Let the environment Hilbert space $\H$ be given, such
%that $\H_\mA\subset\H$ and $\H_M\subset\H$. Given an operator $\mA:\d(A)\to
%\H_{\mathbf{A}}$, that is positive definite in
%the subspace $\H_{\mathbf{A}}\subset\H$, we use
%$$
%P_{\H_{\mathbf{A}}}\mA P_{\H_{\mathbf{A}}}
%$$
%to denote its extension by zero to the whole of $\H$.
%Analogously, we write $P_{\H_M}\mM P_{\H_M}$ for the bounded positive definite operator $M$.

A bounded operator $H:\H\to\H$ is a
Hilbert--Schmidt operator if $H^*H$ is trace class and then, cf. \cite[Ch. X.1.3]{Kato76},
\begin{equation}\label{eq:hs_norm}
\tripleb H \tripleb_{HS}:=\textsf{ Tr} \sqrt{H^*H}.
\end{equation}
%when there exists an orthonormal
%basis\footnote{Such a basis does not have to be countable, cf.
%\cite{\Weidmann_I}.} $(\psi_\alpha)_\alpha$ of $\H$ such that
%$\sum_\alpha\|H\psi_\alpha\|^2<\infty$, cf. \cite{\Weidmann_I}. Futhermore, operator $H$ is
%Hilbert--Schmidt if and only if
%\begin{equation}\label{eq:hs_norm}
%\tripleb H \tripleb_{HS}:=\Big(\sum_i s_i(H)^2\Big)^{1/2}=
%\Big(\sum_\alpha\|H\psi_\alpha\|^2\Big)^{1/2} <\infty,
%\end{equation}
%where $s_i(H)$ are all the \textit{singular values} of the
%operator $H$. The norm $\tripleb \cdot \tripleb_{HS}$
%is called the Hilbert--Schmidt norm. Equivalently, operator $H:\H_1\to\H_2$ is a Hilbert--Schmidt operator
%when there exist bases $(\psi_\alpha)_\alpha$ and $(\phi_\beta)_\beta$ in $\H_1$ and in $\H_2$, respectively,
%such that $\sum_{\alpha, \beta}|(\phi_\beta, H\psi_\alpha)|^2<\infty$. Then it holds
%$$
%\tripleb H\tripleb_{HS}=\Big(\sum_i s_i(H)^2\Big)^{1/2}=\Big(\sum_{\alpha, \beta}|(\phi_\beta, H\psi_\alpha)|^2\Big)^{1/2}.
%$$

Let $\mathbf{A}$ and $\mM$ be positive definite operators in $\H_{\mathbf{A}}\subset\H$
and $\H_\mM\subset\H$, respectively.
We will analyze the weakly formulated Sylvester equation under the
assumption that $\tripleb F\tripleb_{HS}<\infty$ and
\begin{equation}\label{eq:prvi_korak}
\textsf{gap}(\sigma(\mM),\sigma(\mH)):=\inf_{\substack{\mu\in\sigma(\mM),\\
\lambda\in\sigma(\mA)} }
\frac{|\mu-\lambda|}{\sqrt{\mu\lambda}}>0.
\end{equation}

To prove our result, we will need a basic result on the spectral representation of
selfadjoint operators, see \cite[Satz 8.17]{\Weidmann_I}.
\begin{theorem}[Spectral representation]\label{thm:spectral_rep}
For every selfadjoint operator $\mH$ in a sparable Hilbert space $\H$ there
exists a $\sigma$-finite  measure space $(\mathcal{M}, \mu)$, a $\mu$-mea\-sur\-able
function $h:\mathcal{M}\to\R$ and a unitary operator $V:\H\to L^2(\mathcal{M}, \mu)$
such that
$$
\mH=V^{-1} \widetilde{\mH} V.
$$
Here $\widetilde{\mH}: L^2(\mathcal{M}, \mu)\to L^2(\mathcal{M}, \mu)$ is the
multiplication operator which is defined by the function $h$.
\end{theorem}
We will also need the following theorem on the integral
representation of Hilbert--Schmidt operators. For the proof see \cite[Satz 3.19]{\Weidmann_I}.

\begin{theorem}\label{thm:integral_rep}
A bounded operator $T:L^2(\mathcal{M}_1, \mu)\to L^2(\mathcal{M}_2, \nu)$ is a Hilbert--Schmidt
operator if and only if there exists a function $t\in L^2(\mathcal{M}_1\times\mathcal{M}_2, \mu\times\nu)$
such that
$$
(Tg)(y)=\int_{\mathcal{M}_1}t(x,y)g(x) d\mu\qquad \nu
\text{{\rm -almost everywhere, }}\quad g\in L^2(\mathcal{M}_1, \mu).
$$
Furthermore, we have
$$
\tripleb T \tripleb_{HS}=\|t\|_{L^2(\mathcal{M}_1\times\mathcal{M}_2, \mu\times\nu)}.
$$
\end{theorem}

We now prove a ``Hilbert--Schmidt'' version of Theorem \ref{prvo2:t_slabisylv}. We will
assume that $\tripleb F\tripleb_{HS}<\infty$ and that $\H$ be separable. On the other hand,
the spectra of $\mA$ and $\mM$ may be arbitrarily interlaced.

\begin{theorem}\label{cor:unitary_invariant_2}
Let $\mathbf{A}$ and $\mM$ be positive definite operators in $\H_{\mathbf{A}}$
and $\H_{\mM}$, respectively and let $F:\H_{\mM}\to\H_{\mA}$ be a bounded operator.
Assume further that
$\tripleb F\tripleb_{HS}<\infty$ and $\text{{\rm \textsf{gap}}}(\sigma(\mM),\sigma(\mH))>0$
then there exists a unique Hilbert--Schmidt operator $T$
such that
\begin{equation}\label{prvo2:e_slabaHS1}
\sk{\mathbf{A}^{1/2} v, T\mM^{-1/2}u}-\sk{v,\mathbf{A}^{-1/2}T\mM^{1/2}u}=\sk{v, F u}
\end{equation}
and
\begin{equation}\label{eq:ocjena}
\tripleb T \tripleb_{HS} \leq
\frac{\tripleb F\tripleb_{HS}}{\text{{\rm \textsf{gap}}}(\sigma(\mM),\sigma(\mH))}.
\end{equation}
\end{theorem}
\begin{proof}
The uniqueness of the bounded solution of the equation (\ref{prvo2:e_slabaHS1}) follows
by a double cut-off argument analogous to the one used in (\ref{drugo:standardsylv})--(\ref{prvo:e_E_nW}).
We leave out the details.

%Obviously, the assumption $\textsf{gap}(\sigma(\mM), \sigma(\mH))>0$ imples that
%\begin{equation}\label{eq:prvi_korak}
%\sup_{\substack{\mu\in\sigma(\mM),\\
%\lambda\in\sigma(\mA)} }
%\frac{\sqrt{\mu\lambda}}{|\mu-\lambda|}=\frac{1}{\textsf{gap}(\sigma(\mM), \sigma(\mH))}<\infty.
%\end{equation}
By Theorem \ref{thm:integral_rep} there exist measure spaces $(\mathcal{M}_{\mM}, \mu)$
and $(\mathcal{M}_{\mA}, \mu)$, measurable functions $m:\mathcal{M}_{\mM}\to\R$ and
$a:\mathcal{M}_{\mA}\to\R$ and unitary operators $U:\H\to L^2(\mathcal{M}_{\mM}, \mu)$
and $V:\H\to L^2(\mathcal{M}_{\mA}, \mu)$ such that
\begin{align*}
\mA&=V^{-1} \widetilde{\mA} V\\
\mM&=U^{-1} \widetilde{\mM} U.
\end{align*}
Here we have taken $\widetilde{\mA}$ and $\widetilde{\mM}$ to be the multiplication operators which
were defined by the functions $a$ and $m$ respectively.
Since $\tripleb F \tripleb_{HS}<\infty$, the operator
$VFU:L^2(\mathcal{M}_{\mM}, \mu)\to L^2(\mathcal{M}_{\mA}, \mu)$ is obviously a
Hilbert--Schmidt operator and $\tripleb VFU\tripleb_{HS}=\tripleb F\tripleb_{HS}$. We can therefore
assume, without loosing generality, that we work with $\H_{\mM}=L^2(\mathcal{M}_{\mM}, \mu)$,
$\H_{\mA}=L^2(\mathcal{M}_{\mH}, \nu)$ and that $\mA=\widetilde{\mA}$, $\mM=\widetilde{\mM}$
and $F=VFU$.

Theorem \ref{thm:integral_rep} implies that there exists a function
$f\in L^2(\mathcal{M}_{\mM}\times\mathcal{M}_{\mA}, \mu\times\nu)$ such that
$$
(F g)(y)=\int_{\mathcal{M}_{\mM}}f(x,y)g(x) d\mu\qquad \nu
\text{{\rm -almost everywhere, }}\quad g\in L^2(\mathcal{M}_{\mM}, \mu).
$$
Set
\begin{equation}\label{eq:rewritten}
t(x,y)=\frac{f(x,y)}{\frac{a(y)^{1/2}}{m(x)^{1/2}}-\frac{m(x)^{1/2}}{a(y)^{1/2}}},
\qquad \mu\times\nu\text{{\rm -almost everywhere}}.
\end{equation}
Relation (\ref{eq:prvi_korak}) and the positive definiteness of
$\mA$ and $\mM$ imply that
$$
\|\frac{a(\cdot)^{1/2}m(\cdot\cdot)^{1/2}}{a(\cdot)-m(\cdot\cdot)}\|_
{L^\infty(\mathcal{M}_{\mM}\times\mathcal{M}_{\mA}, \mu\times\nu)}\leq
\frac{1}{\textsf{gap}(\sigma(\mM), \sigma(\mA))}
$$
thus $t\in L^2(\mathcal{M}_{\mM}\times\mathcal{M}_{\mA}, \mu\times\nu)$ and
\begin{equation}\label{eq:ocjena_1HHS}
\|t\|_{L^2(\mathcal{M}_{\mM}\times\mathcal{M}_{\mA}, \mu\times\nu)}\leq
\frac{1}{\textsf{gap}(\sigma(\mM), \sigma(\mA))}\|f\|_{L^2(\mathcal{M}_{\mM}\times\mathcal{M}_{\mA}, \mu\times\nu)}.
\end{equation}
Now (\ref{eq:rewritten}) can be rewritten as
\begin{equation}\label{eq:funkcijskaHS}
a(y)^{1/2}t(x,y)m(x)^{-1/2}-a(y)^{-1/2}t(x,y)m(x)^{1/2}=f(x,y)
\end{equation}
The kernel $t$ defines a Hilbert--Schmidt operator $T$ with
$$
(v, Tu)=\int \overline{v(y)}t(x,y)u(x)d\mu ~d\nu.
$$
By taking integrals for $v\in\d(\mA^{1/2})$ and $u\in\d(\mM^{1/2})$ we
establish that the equation (\ref{eq:funkcijskaHS}) is equivalent to
(\ref{prvo2:e_slabaHS1}) and the estimate (\ref{eq:ocjena_1HHS})  implies (\ref{eq:ocjena}).
\end{proof}

\subsection{Estimates by other unitary invariant operator norms}
Let $\mathcal{L}(\H)$ be the algebra of all bounded operators
on the Hilbert space $\H$. We will consider \textit{symmetric norms}
 $ \tripleb \cdot \tripleb $ on a subspace $\SS$ of $\lp(\H)$.
 To say that the norm is symmetric on $\SS\subset\lp(\H)$ means that, beside the
usual properties of any norm, it additionally satisfies:
\begin{description}
  \item[(i)] If $B\in\SS$, $A,C\in\lp(\H)$
  then $ABC\in\SS$ and
  $$
   \tripleb ABC \tripleb \leq\|A\| \tripleb B \tripleb \|C\|.
  $$
  \item[(ii)] If $A$ has rank $1$ then $ \tripleb A \tripleb =\|A\|$, where $\|\cdot\|$
  always denotes the standard operator norm on $\lp(\H)$.
  \item[(iii)] If $A\in\SS$ and $U,V$ are unitary on $\H$, then
  $UAV\in\SS$ and $ \tripleb UAV \tripleb = \tripleb A \tripleb $.
  \item[(iv)] $\SS$ is complete under the norm $ \tripleb \cdot \tripleb $.
%  \item[(v)]
%  Additionally, assume $A_n\in\SS$
%  and $\sup \tripleb A_n \tripleb <\infty$ and $A=\wlim_n A_n$ then $A\in\SS$  and
%  $$
%   \tripleb A \tripleb \leq\sup \tripleb A_n \tripleb .
%  $$
\end{description}
The subspace $\SS$ is
defined as a $ \tripleb \cdot \tripleb $--closure of the set of all degenerate operators
in $\lp(\H)$. Such $\SS$ is an ideal in the algebra $\lp(\H)$, cf. \cite{GohbergKrein,SimonTrace}.
Symmetric norms were used in \cite{BathDavMcInt85} in the context of subspace estimates.
If we assume, additionally to the assumptions of Theorem \ref{prvo2:t_potprostori} that
$\tripleb F\tripleb<\infty$ then there exists a unique bounded solution $T$
of the weak Sylvester equation and
$$
\tripleb T\tripleb \leq\frac{\sqrt{D\normm{M}}}{D-\normm{M}}\tripleb F\tripleb.
$$
We now prove this fact.
\begin{theorem}\label{t:usesDK}
Let $\mA$ and $M$ be the selfadjoint operators which satisfy the assumptions
of Theorem \ref{prvo2:t_slabisylv} and let the symmetric norm $\tripleb \cdot\tripleb$
have the property
\begin{quote}
($\mathcal{P}$)
If $\sup \tripleb A_n \tripleb <\infty$ and $A=\wlim_n A_n$ then $A\in\SS$  and
  $$
   \tripleb A \tripleb \leq\sup \tripleb A_n \tripleb .
  $$
\end{quote}
If $\tripleb F\tripleb<\infty$ then there exists a unique bounded operator $T$
such that
$$
\sk{\mathbf{A}^{1/2} v, TM^{-1/2}u}-\sk{v,\mathbf{A}^{-1/2}TM^{1/2}u}=\sk{v, F u}
$$
and
$$
\tripleb T\tripleb \leq\frac{\sqrt{D\normm{M}}}{D-\normm{M}}\tripleb F\tripleb.
$$
\end{theorem}
\begin{proof}
The proof follows by a cut-off argument. We (re)use the construction which was
used in (\ref{drugo:standardsylv}).
Let $f_n(\mathbf{A}^{1/2})$ and $E_n$ be as in (\ref{drugo:standardsylv}). The equation
\begin{equation}\label{eq:standardSyl}
\sk{f_n(\mathbf{A}^{1/2})v, T_nM^{-1/2}u}-\sk{f_n(\mathbf{A}^{1/2})^{-1}v,T_nM^{1/2}u}=\sk{v, E_nF u}
\end{equation}
can now be written as the standard Sylvester equation
$$
f_n(\mathbf{A}^{1/2})^2 T_n-T_nM=f_n(\mathbf{A}^{1/2})E_nFM^{1/2}
$$
which has the unique bounded solution $T_n:\H_M\to\ra(E_n)$ and $\tripleb T_n\tripleb<\infty$
(this follows from \cite[Theorem 5.2]{DavisKahan70}).
The operator $E_nT$ is bounded and satisfies the equation (\ref{eq:standardSyl}) therefore $T_n=E_nT$.
Here we have tacitly assumed $\lp(\H_{\mA})\subset\lp(\H)$.
Furthermore,
\begin{equation}\label{eq:dav_kah_li}
\mathbf{A}^{1/2} E_nTM^{-1/2}-\mathbf{A}^{-1/2}E_nTM^{1/2}=E_nF.
\end{equation}
 We compute,
using the property \textbf{(i)},
\begin{align*}
\tripleb \mathbf{A}^{-1/2}E_nTM^{1/2}\tripleb&\leq\frac{\|M^{1/2}\|}{\sqrt{\|\mA^{-1}\|^{-1}}}\tripleb E_n T\tripleb\\
\tripleb \mathbf{A}^{1/2}E_nTM^{1/2}\tripleb&\geq\frac{\sqrt{\|\mA^{-1}\|^{-1}}}
{\|M^{1/2}\|}\tripleb E_n T\tripleb.
\end{align*}
From these estimates and (\ref{eq:dav_kah_li}) we obtain the uniform upper bound
\begin{equation}\label{eq:uniformBound}
\tripleb E_nT\tripleb \leq\frac{\sqrt{D\normm{M}}}{D-\normm{M}}\tripleb E_nF\tripleb
\leq\frac{\sqrt{D\normm{M}}}{D-\normm{M}}\tripleb F\tripleb.
\end{equation}
Since $E_nT\to T$ in the strong operator topology, Property ($\mathcal{P}$) and the uniform bound
(\ref{eq:uniformBound}) imply $\tripleb T \tripleb<\infty$ and the desired estimate follows.
\end{proof}

\section{Perturbations of spectral subspaces}\label{prvo:s_potprostori}
When comparing two spectral subspaces of operators $\mH$ and $\mM$, which
satisfy (\ref{eq:prva_ocjena}), we have to make an additional assumption on
the location of the spectra. Namely we assume that there exist $D_1<D_2$ such that
the interval $[D_1, D_2]\subset\R$ is contained in the resolvent sets of both $\mH$
and $\mM$.
Let $Q=E_\mH(D_1)$ and $P=E_\mM(D_1)$. We want to estimate the norm of $P-Q$.
The following description of a relation between a pair of orthogonal
projections in a Hilbert space will be sufficient for our considerations.
For the proof see \cite{Kato76}.
\begin{theorem}[Kato]\label{prvo:KatoPQ}
Let $P$ and $Q$ be two orthogonal projections such that
$$
\|P(\I-Q)\|<1.
$$
Then we have the following alternative. Either
\begin{enumerate}
    \item $\ra(P)$ and $\ra(Q)$ are isomorphic and
    $$
    \|P(\I-Q)\|=\|Q(\I-P)\|=\|P-Q\|\qquad \text{or}
    $$
    \item $\ra(P)$ is isomorphic to true subspace of $\ra(Q)$
    and
    $$
    \|Q(\I-P)\|=\|P-Q\|=1.
    $$
\end{enumerate}
\end{theorem}
To ease the presentation set $P_\perp=\I-P$ and $Q_\perp=\I-Q$. First, let us consider the
case when $h$ and $m$ are positive definite. With the
help of Proposition \ref{lema:sim} we shall later
reduce the nonnegative definite case to the positive definite
one.

We define the operators
\begin{equation}\label{eq:umetnta}
\mA=Q_\perp H Q_\perp, \quad H=Q \mH Q, \quad M=P\mM P\quad \text{and}\quad \mW=P_\perp\mM P_\perp.
\end{equation}
We shall not notationally distinguish the operators $\mA$, $M$, $\mW$ and $H$ from
their restrictions to the complement of their respective kernels.
Obviously,
$$
\mH=H+\mA,\qquad \mM=M+\mW
$$
and we compute, for $S$ from (\ref{Tgeneral1}),
\begin{align}
\nonumber Q_\perp S P&= (\mH^{1/2}Q_\perp P \mM^{-1/2}- \mH^{-1/2}Q_\perp P\mM^{1/2})P\\
\nonumber &=\mA^{1/2} Q_\perp P M^{-1/2} -\mA^{-1/2} Q_\perp PM^{1/2}\\
&=\mA^{1/2} T M^{-1/2} -\mA^{-1/2} T M^{1/2}.\label{eq:operator_version}
\end{align}
Here we have defined $T=Q_\perp P$. If we assume that $\textsf{dim}(Q)=\textsf{dim}(P)<\infty$
then Theorem \ref{prvo:KatoPQ} yields
$$
\|P-Q\|=\|T\|.
$$
The case when $\textsf{dim}(Q)=\textsf{dim}(P)=\infty$ will follow in a similar fashion.

The operator equation can be written in the following variational form
\begin{align}\label{prvo:e_zbunjujuca}
(\mA^{1/2} v, TM^{-1/2} u)&-(\mA^{- 1/2} v, TM^{1/2} u)=(v, Su),\\
\nonumber & v\in\d(\mA^{1/2}),\quad u\in \ra(P),
\end{align}
which we have called the weakly formulated Sylvester equation.

\begin{theorem}\label{prvo2:t_potprostori}
Let the positive definite forms $m$ and $h$ be given such that (\ref{eq:prva_ocjena}) holds.
Let there exist $D_1<D_2$ such that
the interval $[D_1, D_2]\subset\R$ be contained in the resolvent sets of both $\mH$
and $\mM$.
Set $Q=E_\mH(D_1)$, $P=E_\mM(D_1)$ and assume $\eta<(D_2-D_1)(D_2 D_1)^{-1/2}$ then
\begin{equation}\label{prvo2:e_pazi_ sad}
\norm{P-Q}\leq\frac{\sqrt{D_2D_1}}{D_2-D_1}~\eta~.
\end{equation}
\end{theorem}
\begin{proof}
$T=Q_\perp P$ is the unique solution of the equation (\ref{prvo:e_zbunjujuca}).
Theorem \ref{cor:unitary_invariant} implies
\begin{eqnarray}
\nonumber \|T\|&\leq&\frac{\eta}{2}\sqrt{\frac{D_2\lambda_n(\mM)}{(D_2-d)(d-\lambda_n(\mM))}}.
\label{prvo2:e_optimald}
\end{eqnarray}
for any $\lambda_n(\mM)< d < D_2$. The optimal $d$ equals $\frac{(D_2+\lambda_n(\mM))}{2}$ and
since $\|M\|<D_1$ we conclude
$$
\norm{Q_\perp P}\leq \eta~\frac{\sqrt{D_2\lambda_n(\mM)}}{D_2-\lambda_n(\mM)}
\leq\eta~\frac{\sqrt{D_2D_1}}{D_2-D_1}<1.
$$
Analogous argumentation for $T=P_\perp Q$, with the roles of $\mH$ and $\mM$ in (\ref{prvo:e_zbunjujuca})
being interchanged, yields the inequality
$$
\norm{P_\perp Q}\leq \eta~\frac{\sqrt{D_2 D_1}}{D_2-D_1}<1.
$$
Theorem (\ref{prvo:KatoPQ}) now implies that
$$
\|Q_\perp P\|=\|P_\perp Q\|=\|Q-P\|.
$$
This in turn establishes (\ref{prvo2:e_pazi_ sad}).
\end{proof}
In the case in which $h$ is only nonnegative definite, assumption (\ref{eq:prva_ocjena}) implies that
$\je(\mM)=\je(\mH)$ and $\ra(\mH)=\ra(\mM)$, since $\mH$ and $\mM$ are selfadjoint.
This in turn allows us to conclude
that $\mathcal{N}:=\ra(P)\cap\je(\mH)\subset\ra(Q)$, so
$$
\widetilde{Q}=Q-P_{\mathcal{N}}, \qquad \widetilde{P}=P-P_{\mathcal{N}}
$$
are orthogonal projections and
$$
\|Q-P\|=\|\widetilde{Q}-\widetilde{P}\|.
$$
Since $\ra(\widetilde{P})\subset\ra(\mH)$ and $\ra(\widetilde{Q})\subset\ra(\mH)$ we
can reduce the problem to the positive definite case.
\begin{theorem}\label{prvo1:t_sintheta}
Let the positive definite forms $m$ and $h$ be given such that (\ref{eq:prva_ocjena}) holds.
Let there exist $0<L_1<L_2<D_1<D_2$ such that
the intervals $[L_1, L_2]\subset\R$ and $[D_1, D_2]\subset\R$ be
contained in the resolvent sets of both $\mH$
and $\mM$.
Set $Q=E_\mH[L_1, L_2]$, $P=E_\mM[D_1, D_2]$ and assume $$\Big[\frac{\sqrt{D_2D_1}}{D_2-D_1}+\frac{\sqrt{L_2L_1}}{L_2-L_1}\Big]~\eta<1$$ then
\begin{equation}\label{prvo2:e_pazi_ sad2}
\norm{P-Q}\leq\Big[\frac{\sqrt{D_2D_1}}{D_2-D_1}+\frac{\sqrt{L_2L_1}}{L_2-L_1}\Big]~\eta~.
\end{equation}
\end{theorem}
\begin{proof}
The assumption $L_1>0$ implies that we may assume, without losing any generality,
that we have the positive definite forms $m$ and $h$.
Theorem \ref{cor:unitary_invariant_2} and the same argument
as in Theorem \ref{prvo2:t_potprostori} implies
$$
\|P_\perp Q\|=\|Q_\perp P\|=\|P-Q\|.
$$
This in turn allows us to conclude that
\begin{eqnarray}
\norm{P-Q}\leq\Big[\frac{\sqrt{D_2D_1}}{D_2-D_1}+\frac{\sqrt{L_2L_1}}{L_2-L_1}\Big]~\eta.\label{prvo2:e_dupli}
\end{eqnarray}
\end{proof}
Numerical experiments with the Sturm--Liouville eigenvalue problem, which
were performed in \cite{GruVes02}, illustrated that in some situations
the results of Theorems \ref{prvo2:t_potprostori}
and \ref{prvo1:t_sintheta} yield considerably sharper estimates of the perturbations of the
spectral subspaces than the results of \cite{BathDavMcInt85,DavisKahan70}. We now show
that our theorems also apply in situations in which the theory from
\cite{BathDavMcInt85,DavisKahan70} does not.
\begin{example}\label{example:measur}
Take \(\mH\), \(\mM\) as selfadjoint realizations of the differential
operators
\[
-\frac{\partial}{\partial x}\alpha(x)\frac{\partial}{\partial x},\quad
-\frac{\partial}{\partial x}\beta(x)\frac{\partial}{\partial x},
\]
respectively, in the Hilbert space \(\H=L^2(I)\), \(I\) a (finite or infinite)
interval with, say, Dirichlet  boundary conditions and non-negative
bounded measurable functions $\alpha(x)$, $\beta(x)$ which satisfy
\[
|\beta(x) - \alpha(x)| \leq \eta\sqrt{\beta(x)\alpha(x)}.
\]
Now, the form
$$
\delta(u,v)=h(u, v)-m(u, v)
$$
is not---in general---representable by a bounded operator. This rules out an application
of the subspace perturbation theorems from \cite{BathDavMcInt85,DavisKahan70}. On the other hand
both of our Theorems \ref{prvo2:t_potprostori} and \ref{prvo1:t_sintheta} apply and
yield the corresponding estimates, e.g.
$$
\norm{E_\alpha(D_1)-E_\beta(D_1)}\leq\frac{\sqrt{D_2D_1}}{D_2-D_1}~\eta~,
$$
when we know that $[D_1, D_2]$ is contained in the resolvent sets of both $\mH$ and $\mM$.
\end{example}

Theorem \ref{cor:unitary_invariant_2} can also be applied to
yield a Hilbert--Schmidt version of Theorems \ref{prvo2:t_potprostori} and \ref{prvo1:t_sintheta}.

\begin{theorem}\label{prvo1:t_sinthetaHS}
Let the positive definite forms $m$ and $h$ be given such that (\ref{eq:prva_ocjena}) holds.
Assume $P$ and $Q$ are projections which commute with the operators $\mH$ and $\mM$
respectively and let $\mA$, $M$, $\mathbf{W}$, $H$ as in (\ref{eq:umetnta}).
If $\tripleb Q_\perp SP \tripleb_{{\rm HS}}<\infty$,
$\tripleb P_\perp S^* Q\tripleb_{{\rm HS}}<\infty$ and
$$
\text{{\rm \textsf{gap}}}(\sigma(\mA), \sigma(M)),\qquad \text{{\rm \textsf{gap}}}(\sigma(\mathbf{W}), \sigma(H)),
$$
are both positive.
Then $Q_\perp P$, $P_\perp Q$ and $P-Q$ are Hilbert--Schmidt operators and
\begin{align}\label{eq:prvaHsQP}
\tripleb Q_\perp P\tripleb_{{\rm HS}}^2&\leq\frac{\tripleb Q_\perp SP \tripleb_{{\rm HS}}^2}
{\text{{\rm \textsf{gap}}}(\sigma(\mA), \sigma(M))^2}\\
\tripleb P_\perp Q\tripleb_{{\rm HS}}^2&\leq\frac{\tripleb Q S P_\perp\tripleb^2_{{\rm HS}}}{\text{{\rm \textsf{gap}}}(\sigma(\mathbf{W}), \sigma(H))^2}\\
\tripleb P-Q\tripleb_{{\rm HS}}^2&\leq\frac{\tripleb Q_\perp SP \tripleb_{{\rm HS}}^2}
{\text{{\rm \textsf{gap}}}(\sigma(\mA), \sigma(M))^2}+
\frac{\tripleb Q S P_\perp\tripleb^2_{{\rm HS}}}{\text{{\rm \textsf{gap}}}(\sigma(\mathbf{W}), \sigma(H))^2}~.
\label{prvo2:e_pazi_ sad2HS}
\end{align}
\end{theorem}
\begin{proof}
by construction (\ref{eq:umetnta}) the operator $T=T_1=Q_\perp$ satisfies
the Sylvester equation (\ref{prvo:e_zbunjujuca}), which in this setting has the form, cf. (\ref{prvo2:e_slabaHS1}),
\begin{align}\label{prvo:e_zbunjujuca2}
(\mA^{1/2} v, T M^{-1/2} u)&-(\mA^{- 1/2} v, T M^{1/2} u)=(v, Q_\perp S Pu),\\
\nonumber & v\in\d(\mA^{1/2}),\quad u\in \d( M^{1/2}).
\end{align}
On the other hand, the operator
$T=T_2=P_\perp Q$ satisfies the ``dual'' equation, cf. Proposition \ref{lema:sim},
\begin{align}\label{prvo:e_zbunjujuca3}
(\mW^{1/2} v, T H^{-1/2} u)&-(\mW^{- 1/2} v, T H^{1/2} u)=(v, P_\perp S^* Qu),\\
\nonumber & v\in\d(\mW^{1/2}),\quad u\in\d(H^{1/2}).
\end{align}
Now,
$$
(P-Q)^2=Q_\perp P+P_\perp Q,
$$
where by Theorem \ref{cor:unitary_invariant_2} both $Q_\perp P$ and $P_\perp Q$
are Hilbert--Schmidt\footnote{To prove this equality
one can use the singular value analysis from \cite{DavisKahan70}. Alternatively,
one could use the property ($\mathcal{P}$) from Theorem \ref{t:usesDK}.} and
\begin{align*}
\tripleb P-Q\tripleb_{HS}^2&=\textsf{Tr}(Q_\perp P+P_\perp Q)=
\textsf{Tr}(P Q_\perp P)+\textsf{Tr}(Q P_\perp Q)\\
&=\tripleb Q_\perp P\tripleb_{HS}^2+\tripleb P_\perp Q\tripleb_{HS}^2
\end{align*}
Using (\ref{eq:ocjena}), we see that estimates (\ref{eq:prvaHsQP})--(\ref{prvo2:e_pazi_ sad2HS})
hold.
\end{proof}

\begin{corollary}\label{prvo1:c_sinthetaHS}
Let the positive definite forms $m$ and $h$ be given such that (\ref{eq:prva_ocjena}) holds.
If $\tripleb S\tripleb_{{\rm HS}}<\infty$ and the other conditions of Theorem \ref{prvo1:t_sinthetaHS}
hold. Then
\begin{equation}\label{eq:estimateHSC}
\tripleb P-Q\tripleb_{{\rm HS}}^2\leq\frac{\tripleb S \tripleb_{{\rm HS}}^2}
{\min\big\{\text{{\rm \textsf{gap}}}(\sigma(\mA), \sigma(M))^2, \text{{\rm \textsf{gap}}}(\sigma(\mathbf{W}), \sigma(H))^2\big\}}~.
\end{equation}
\end{corollary}
\begin{proof}
Just note that $$\tripleb Q_\perp SP\tripleb_{HS}^2+\tripleb QSP_\perp\tripleb_{HS}^2\leq
\tripleb S\tripleb_{HS}^2.
$$
\end{proof}

\section{Further properties of the operator $S$ --- an application in the numerical analysis}
We will now present an application of Theorem
\ref{cor:unitary_invariant} in numerical analysis. This will also
demonstrate a role played by the new Hilbert--Schmidt norm
estimates.

Assume now that we are given a positive definite operator $\mH$
such that (\ref{eigenvaluesH}) holds. Let $P$ be an orthogonal
projection such that $\ra(P)\subset \q(h)$ and $\textsf{dim}~
\ra(P)=n$. We aim to obtain estimates of
\begin{equation}\label{eq.ciljocj}
\tripleb E_{\mH}(D)-P\tripleb
\end{equation}
for $\tripleb\cdot\tripleb=\|\cdot\|$ and
$\tripleb\cdot\tripleb=\tripleb\cdot\tripleb_{HS}$.

We estimate (\ref{eq.ciljocj}) by an application of Theorem
\ref{cor:unitary_invariant} (equivalently Theorem
\ref{prvo1:t_sinthetaHS}). Theorem \ref{cor:unitary_invariant}
will allow us to improve \cite[Theorem 3.2]{GruVes02} inasmuch as
that we establish estimates for the Hilbert--Schmidt norm  and not
just the spectral norm.

The properties of the main perturbation construction from
\cite{GruVes02}, cf. \cite{Gru03_3,GruPhd}, will be summarized for
reader's convenience.

We start by defining the positive definite form
\begin{equation}\label{eq_p_hp}
h_P(u,v)=h(Pu, Pv)+h(P_\perp u, P_\perp v)
\end{equation}
and a self adjoint operator $\mH_P$ which represents the form $h_P$ in the sense of Kato.
It was shown, see \cite{Gru03_3,GruPhd,GruVes02}
that
\begin{enumerate}
\item the form $h_P$ is positive definite , hence there exists the positive definite operator
$\mH_P$ which represents $h_P$ in the sense of Kato.
\item $\q(h)=\q(h_P)$
\item $\ra(P)$ reduces $\mH_P$.
\item $\mH^{-1}-\mH_P^{-1}$ is a degenerate self adjoint operator.
\item Let $\delta h^P:=h-h_P$ and let $\delta H_s^P$ be the bounded self adjoint operator
such that
\begin{equation}\label{eq.deltahs}
(u,\delta H_s^P v)=\delta h^P(\mH_P^{-1/2}u,\mH_P^{-1/2}v)
\end{equation}
then $\delta H_s^P$ is a degenerate operator and $\textsf{dim} \ra(\delta H_s^P)=2n$.
\item The values
\begin{equation}\label{eq.sing_val}
\eta_i=\max_{\substack{\mathcal{S}\subset\ra(P)\\
\textsf{dim}\mathcal{S}=n-i}}\min\Big\{
\frac{(\psi,\mH^{-1}\psi)-(\psi,\mH^{-1}_P\psi)}{(\psi,\mH^{-1}\psi)}:\psi\in\mathcal{S}\Big\}
\end{equation}
together with their negatives are all non-zero eigenvalues of $\delta H_s^P$. Furthermore,
$\eta_i$ are all the singular values of the operator $\delta H_s^P P$.
\item \begin{equation}\label{eq.sin_theta_gv}
|\delta h^P(\phi,\psi)|\leq\eta_n\sqrt{h_P[\psi]h_P[\phi]}
\end{equation}
\end{enumerate}

The estimates from  \cite[Theorem 3.2]{GruVes02} only use information which is contained in $\eta_n$.
New theory allows us to take advantage of other $\eta_i$.

\begin{proposition}\label{prop.simnorm}
Let $P$ and $h_P$ be as in (\ref{eq_p_hp}) and let
$$
S=\mH^{1/2}\mH^{-1/2}_P-\overline{\mH^{-1/2}\mH^{1/2}_P}
$$
then
\begin{align*}
\tripleb SQ\tripleb&\leq\frac{\tripleb\delta H_s^P Q\tripleb}{\sqrt{1-\eta_n}}\\
\tripleb S\tripleb&\leq\frac{\tripleb\delta H_s^P \tripleb}{\sqrt{1-\eta_n}}.
\end{align*}
Here $\delta H_s^P$ is the degenerate operator from (\ref{eq.deltahs}), $Q$ is any projection,
$\eta_n$ is given by (\ref{eq.sing_val}) and $\tripleb\cdot\tripleb$ is any
unitary invariant norm.
\end{proposition}
\begin{proof}
\begin{align*}
(\psi, S Q\phi)&=\delta h^P(\mH^{-1/2}\psi, \mH_P^{-1/2}  Q \phi)=
\delta h^P(\mH^{-1/2}_P\big(\mH^{1/2}_P\mH^{-1/2}\big)\psi, \mH^{-1/2}_P Q \phi)\\
&=(\psi,\big(\mH^{1/2}_P\mH^{-1/2}\big)^*\delta H_s^P Q \phi),\qquad \phi,\psi\in\H
\end{align*}
(\ref{eq.sin_theta_gv}) and (\ref{prvo:e_ocjenaZaUpotr}) imply $\|\mH^{1/2}_P\mH^{-1/2}\|\leq 1/\sqrt{1-\eta_n}$.
Property \textbf{(i)} of the symmetric norm $\tripleb\cdot\tripleb$ and the fact that $\delta H_s^P Q$ is a
degenerate operator allow us to complete the proof.
\end{proof}
This proposition leads to an
improved version of \cite[Theorem 3.2]{GruVes02}. Observe that $\tripleb\delta H_s^P\tripleb$ depends only
on $\eta_i$ from (\ref{eq.sing_val}).
\begin{theorem}\label{prvo2:t_potprostori_symm}
Let $h$ be as in (\ref{eigenvaluesH}) and let $P$ and $h_P$ be as in (\ref{eq_p_hp})
and $\eta_i$ as in (\ref{eq.sing_val}). Set
$$
D_P:=\max_{\psi\in\ra(P)}\frac{h[\psi]}{~\|\psi\|^2}
$$
and assume $\eta_n(1-\eta_n)^{-1}<(D-D_P)(D+ D_P)$ then
\begin{equation}\label{prvo2:e_pazi_ sad_2}
\tripleb E_{\mH}(\lambda_n)P_\perp \tripleb\leq\frac{\sqrt{DD_P}}{D-D_P}~\frac{\tripleb\delta H_s^P P\tripleb}{\sqrt{1-\eta_n}}~.
\end{equation}
Here $\tripleb\cdot\tripleb$ is any unitary invariant norm which has Property ($\mathcal{P}$).
\end{theorem}
\begin{proof}
Set $T=(E_\mH(\lambda_n))_\perp P$ and apply Theorem \ref{t:usesDK} to
estimate the norm \hfill\\$\tripleb(E_\mH(\lambda_n))_\perp P\tripleb$. Proposition \ref{prop.simnorm} now implies (\ref{prvo2:e_pazi_ sad_2}), cf.
Corollary \ref{prvo1:c_sinthetaHS}, \cite[Corollary 3.1]{DavisKahan70} and \cite[Proposition 6.1]{DavisKahan70}.
\end{proof}

Assume $\tripleb\cdot\tripleb=\tripleb\cdot\tripleb_{HS}$, then Theorem \ref{prvo2:t_potprostori_symm}
yields the estimate
\begin{equation}\label{prvo2:e_pazi_ sad_HS}
\tripleb (E_{\mH}(\lambda_n))_\perp P\tripleb_{HS}\leq\frac{\sqrt{DD_P}}{D-D_P}~
\frac{\sqrt{\eta_1^2+\cdots+\eta_n^2}}{\sqrt{1-\eta_n}}~.
\end{equation}
\begin{remark}
Note that $\|(E_{\mH}(\lambda_n))_\perp P\|=
\|E_{\mH}(\lambda_n)- P\|$, cf. Theorem \ref{prvo:KatoPQ}. A similar relation holds for a general unitary invariant norm since
according to \cite[Corollary 3.1]{DavisKahan70} and \cite[Section 2]{DavisKahan70}
we have $\tripleb (E_{\mH}(\lambda_n))_\perp P\tripleb=\tripleb P_\perp E_{\mH}(\lambda_n)\tripleb$ and
\begin{equation}\label{eq:samas}
\tripleb E_{\mH}(\lambda_n)-P\tripleb=\tripleb (E_{\mH}(\lambda_n))_\perp P+ P_\perp E_{\mH}(\lambda_n)\tripleb .
\end{equation}
An estimate of (\ref{eq:samas}) is obtained by
a combination of Proposition \ref{prop.simnorm} and available (depending on an application)
information on the separation of the involved spectra, cf. Corollary \ref{prvo1:c_sinthetaHS}.
We have not specified a general estimate on $\tripleb E_{\mH}(\lambda_n)-P \tripleb$ since we
consider such estimates to be highly application dependent and we would not like to prejudice
their form.
\end{remark}
\begin{table}[h]
\begin{center}
\begin{tabular}{|l||c|ccccc|}\hline
%&&&&&&&\cr
$N$ &5&6&7&8&9&10\cr
\hline
&&&&&&\cr
$\displaystyle \tripleb (E_{\mH}(\lambda_2))_\perp P\tripleb_{HS}$
&4.4e-3 &  2.0e-3 & 1.1e-3 & 6.0e-4  & 3.7e-4 & 2.4e-4 \\
&&&&&&\cr\hline
&&&&&&\cr
 $\displaystyle \frac{\sqrt{\lambda_3D_{P_N}}}{\lambda_3-D_{P_N}}
 \frac{\sqrt{\eta_1^2+\eta^2_2}}{\sqrt{1-\eta_2}}$
& 2.2e-2 &1.0e-2& 5.3e-3 & 3.3e-3 & 2.2e-3 & 1.5e-3\\
 &&&&&&\cr\hline
 &&&&&&\cr
 $\displaystyle\frac{\sqrt{s_1(R_2^N)+s_2(R_2^N)}}{\lambda_3-D_{P_N}}$
 &2.0e-2&  1.4e-2&   9.6e-3&  7.2e-3&  5.5e-3 &4.4e-3\\
&&&&&&\cr
\hline\end{tabular}
\end{center}
\caption{Error estimate from Theorem \ref{prvo2:t_potprostori_symm} and the true error}
\label{tab:prva}
\end{table}

We will now evaluate (\ref{prvo2:e_pazi_ sad_HS}) on the example from \cite[Section 4]{GruVes02}.
There we have considered the positive definite operator $\mH$ which is defined by
the symmetric form
\begin{align*}
h(u,v)&=\int^{2\pi}_0\big(u'\overline{v'}-\alpha u\overline{v} \big)~dt
\label{e_slabiproblem1}\\
\nonumber u, v &\in\{f:f,f'\in L^2[0,2\pi], {\rm e}^{{\rm i}
\theta}f(0)=f(2\pi)\}=\d(h).
\end{align*}
The eigenvalues and eigenvectors of the operator $\mH$ are
\begin{eqnarray*}
\omega_{\pm k}&=&\left(\pm k+\frac{\theta}{2\pi}\right)^2 -\alpha,\quad
z_{\pm k}(t)=\eexp^{-\imag \left(\pm k+\frac{\theta}{2\pi}\right)
t},
\quad k\in\N\\
\omega_{0}&=&\left(\frac{\theta}{2\pi}\right)^2-\alpha,\qquad\quad
z_{0}(t)=\eexp^{- \imag \frac{\theta}{2\pi} t}.
\end{eqnarray*}
In standard notation we have
\begin{align*}
\lambda_1(\mH)&=\omega_0,\;\;\lambda_2(\mH)=\omega_{-1},\;\;
\lambda_3(\mH)=\omega_{1},\\
u_1&=z_{0},\;\;u_2=z_{-1},\;\;u_3=z_{1}.
\end{align*}
For numerical experiments we chose $\theta=\pi-10^{-4}$ and
$\alpha=0.2499$ so that the eigenvalues $\lambda_1$ and
$\lambda_2$ are ``small'' and tightly clustered. As a test space
we chose $\mathcal{Y}^3_N=\textsf{span}\big\{w_1^N, w_2^N\big\}$,
where $w_1^N$ and $w_2^N$ are generated by the \textit{smooth} $N$
point equidistant \textit{cubic} interpolation of the known
eigenfunctions $u_1$ and $u_2$. Take $P_N$ such that
$\ra(P_N)=\mathcal{Y}^3_N$. Since $\mathcal{Y}^3_N\subset\d(\mH)$
both Theorem \ref{prvo2:t_potprostori_symm} and the bounds from
\cite{DavisKahan70} apply. Set $r_\phi =\mH\phi + (\phi,
\mH\phi)\phi$, for $\|\phi\|=1$. Since $w_1^N, w_2^N\in\d(\mH)$ we conclude that
$r_{w_1^N}$ and $r_{w_2^N}$ are bona fide vectors. Set
$$
R_2^N=\begin{bmatrix} (r_{w_1^N}, r_{w_1^N}) & (r_{w_1^N}, r_{w_2^N})\\
(r_{w_2^N}, r_{w_1^N})&(r_{w_2^N}, r_{w_2^N})\end{bmatrix}.
$$
The competing bound from \cite{DavisKahan70} is
\begin{equation}\label{eq:dav.kah.comp}
\tripleb (E_{\mH}(\lambda_2))_\perp P_N\tripleb_{HS}\leq\frac{\sqrt{s_1(R_2)+s_2(R_2)}}{\lambda_3-D_{P_N}}.
\end{equation}
On Table \ref{tab:prva} we have displayed the actual measured
error in the first line, in the second line we display the bound
from (\ref{prvo2:e_pazi_ sad_HS}) and in the third line
Davis--Kahan bound (\ref{eq:dav.kah.comp}). We see that with the
improvement of the approximation the advantage of the bound from
Theorem \ref{prvo2:t_potprostori_symm} over
(\ref{eq:dav.kah.comp}) grows, see Table \ref{tab:prva}. Further
examples, where a numerical advantage of (\ref{prvo2:e_pazi_
sad_2}) over (\ref{eq:dav.kah.comp}) is more stunning, are given
in \cite{GruVes02}. We repeat the results of the numerical
experiments from \cite{GruVes02} on Table \ref{tab:druga}. There
we try to estimate the approximation error in the vector $w^N_1$
in the $\|\cdot\|$-norm by an application of Theorem
\ref{prvo2:t_potprostori}. Otherwise the makeup of Table
\ref{tab:druga} is the same as the makeup of Table \ref{tab:prva}.

\begin{table}[h]
\begin{center}
\begin{tabular}{|l||ccccc|}\hline
%&&&&&&&\cr
$N$ &6&7&8&9&10\cr
\hline
&&&&&\cr
$\displaystyle \| E_{\mH}(\lambda_1)-P_N\|$
&  2.0e-3&  1.1e-3&  6.0e-4&   3.7e-4& 2.4e-4\\
&&&&&\cr\hline
&&&&&\cr
 $\displaystyle \frac{\sqrt{\lambda_2d_{P_N}}}{\lambda_2-d_{P_N}}\frac{\eta_2}{\sqrt{1-\eta_2}}$
& 1.5e0 &  6.2e-1 & 3.5e-1 & 2.2e-1& 1.5e-1\cr
 &&&&&\cr\hline
 &&&&&\cr
 $\displaystyle\frac{\sqrt{s_1(R_2^N)}}{\lambda_2-d_{P_N}}$& 3.6e+2 &  2.1e+2 & 1.5e+2 & 1.1e+2& 8.9e+1\cr
 &&&&&\cr
\hline\end{tabular}
\end{center}
\caption{Approximations for $u_1$ (here we use $d_{P_N}:=\min_{\psi\in\ra(P_N)}\frac{h[\psi]}{~\|\psi\|^2}$)}
\label{tab:druga}
\end{table}

We now present a variation on this example where
(\ref{eq:dav.kah.comp}) does not apply whereas (\ref{prvo2:e_pazi_
sad_2}) still gives useful information. We chose
$\mathcal{Y}^1_N=\textsf{span}\big\{l_1^N, l_2^N\big\}$, where
$l_1^N$ and $l_2^N$ are generated by the $N$ point equidistant
\textit{continuous} \textit{linear} interpolation of $u_1$ and
$u_2$ then $r_{l_1^N}$ and $r_{l_2^N}$ are no longer bona fide
vectors. Subsequently, (\ref{eq:dav.kah.comp}) does not apply any
more but Theorem \ref{prvo2:t_potprostori_symm} is still
applicable. Take now $Q_N$ such that $\ra(Q_N)=\mathcal{Y}^1_N$.
The results are presented on Table \ref{fig.linear}.

\begin{table}[t]
\begin{center}
\begin{tabular}{|l||ccc|}\hline
$N$ &100&120&140\cr\hline
&&&\cr
$\displaystyle \tripleb (E_{\mH}(\lambda_2))_\perp Q_N\tripleb_{HS}$&
5.2024e-005 & 3.6126e-005&  2.6541e-5\\
&&&\cr\hline
&&&\cr
 $\displaystyle
\frac{\sqrt{\lambda_3D_{Q_N}}}{\lambda_3-D_{Q_N}}
\frac{\sqrt{\eta_1^2+\eta_2^2}}{\sqrt{1-\eta_2}}$&
8.7374e-003 & 6.9293e-003 & 5.7302e-003\\
&&& \cr
\hline\end{tabular}
\end{center}
\caption{}
\label{fig.linear}
\end{table}

The performance of the bound (\ref{prvo2:e_pazi_ sad_2}) is influenced by the
quotient
$$
\frac{|D_{P_N}-\lambda_2|}{D_{P_N}}.
$$
$D_{P_N}$ is an approximation\footnote{To be more precise $D_{P_N}$ is
Rayleigh--Ritz approximation to $\lambda_2(\mH)$ from the subspace $\ra(P_N)$.
For more on the Rayleigh--Ritz eigenvalue approximations see \cite{GruVes02}.}
of $\lambda_2$ and in this example we have measured
$$
\frac{|D_{P_N}-\lambda_2|}{D_{P_N}}>0.17,\qquad N=100, 120, 140.
$$
The (under)performance of the bound (\ref{prvo2:e_pazi_ sad_2}) correctly detects this
approximation feature of $\ra(Q_N)$, cf. Table \ref{fig.linear}.

\section{Estimates for perturbations of the square root of a nonnegative operator}
In this section we will show that there are interesting applications of the
equation (\ref{drugo:weak_sylvester}) even when all of the coefficients
$\mA, M$ and $F$ are unbounded.
To demonstrate this we will generalize the known scalar inequality
\footnote{``The relative error in the square root is bounded by
the half relative error in the radicand''.}
\begin{equation}\label{sqrtnumbers}
\frac{|\sqrt{m} - \sqrt{h}|}{\sqrt[4]{mh}} \leq
 \frac{|m - h|}{2\sqrt{mh}},\quad h,m > 0.
\end{equation}
to positive definite Hermitian matrices or, more generally, to
positive, possibly unbounded, operators in an arbitrary Hilbert space.
One of the obtained bounds is
\begin{equation}\label{sqrtmatrices}
\|M^{-1/4}(M^{-1/2} - H^{-1/2})H^{-1/4}\| \leq
\frac{1}{2}\|M^{-1/2}(M - H)H^{-1/2}\|.
\end{equation}
In \cite{Math} a related bound for finite matrices was obtained. It reads
\begin{equation}\label{sqrtMath}
\|H^{-1/4}(M^{-1/2} - H^{-1/2})H^{-1/4}\|\leq
\frac{\eta}{2} + O(\eta^2),
\end{equation}
\[
\eta = \|H^{-1/2}(M - H)H^{-1/2}\|.
\]
This is a more common type of estimate --- the error is measured
by the ``unperturbed operator'' only --- while in our estimate the error
is measured by \(H\) and \(M\) in a symmetric way. The latter type
of estimate is convenient, if both operators \(H\) and \(M\) are
known equally well and we are interested in a possibly sharp bound.
Our bound is obviously as sharp as its scalar pendant. It is also
rigorous, in contrast to (\ref{sqrtMath}) which is only asymptotic.
Moreover, (\ref{sqrtmatrices}) will retain its validity for
fairly general positive selfadjoint operators in a Hilbert space.
The bound (\ref{sqrtmatrices}) is a ``relative bound'' which may be
convenient in computing or measuring
environments (cf. related
bounds obtained for the eigenvalues and eigenvectors of the
Hermitian matrices in \cite{MathVeselic98} and the literature cited there).
Also, this bound is readily
expressed in terms of quadratic forms, which will be convenient
for application with elliptic differential operators as will be
shown below.

The idea of the proof is very simple, especially in the finite
dimensional case which we present first, also in order to accommodate
readers not interested in infinite dimension technicalities.

The basis of our proof is the obvious Sylvester equation
(cf. \cite{Schmitt})
\begin{equation}\label{sylvmatr}
M^{1/2}(M^{1/2} - H^{1/2}) + (M^{1/2} - H^{1/2})H^{1/2} = M - H,
\end{equation}
valid for any Hermitian, positive definite matrices \(H\) and \(M\).
We rewrite this equation in the equivalent form
\begin{equation}\label{oursylvmatr}
M^{1/4}XH^{-1/4} + M^{-1/4}XH^{1/4} = T
\end{equation}
with
\begin{equation}\label{TXmatr}
T = M^{-1/2}(M - H)H^{-1/2},\quad
X = M^{-1/4}(M^{-1/2} - H^{-1/2})H^{-1/4},
\end{equation}
which is immediately verified. This equation
has a unique solution
\begin{equation}\label{solutionmatr}
X = \int_0^\infty e^{-M^{-1/2}t}M^{-1/4}TH^{-1/4}e^{-H^{-1/2}t}dt.
\end{equation}
(just premultiply (\ref{oursylvmatr}) by \(e^{-M^{-1/2}t}M^{-1/4}\),
postmultiply by \(e^{-H^{-1/2}t}H^{-1/4}\), integrate from
\(0\) to \(\infty\) and perform partial integration on its left hand side).
Hence for arbitrary vectors \(\phi\), \(\psi\) we have
\[
|(X\psi,\phi)|^2 \leq \|T\|^2
\left(\int_0^\infty \|e^{-M^{-1/2}t}M^{-1/4}\phi\|
\|e^{-H^{-1/2}t}H^{-1/4}\psi\|dt\right)^2
\]
\[
\leq \|T\|^2 \int_0^\infty \|e^{-M^{-1/2}t}M^{-1/4}\phi\|^2dt
\int_0^\infty \|e^{-H^{-1/2}t}H^{-1/4}\psi\|^2dt
\]
\begin{equation}\label{estimation}
= \frac{\|T\|^2}{4}\|\|\psi\|^2\phi\|^2,
\end{equation}
where we have used the obvious identity
\begin{equation}\label{Cidentity}
\int_0^\infty e^{-2Ct}Cdt = \frac{1}{2}I
\end{equation}
for \(C = H^{-1/2},\ M^{-1/2}\).
Thus, (\ref{sqrtmatrices}) holds true.\\

We now turn to the Hilbert space \(\H\) of arbitrary dimension.
We assume that  \(\mH\) and \(\mM\) are positive selfadjoint operators.
%According to the Kato's monograph \cite{Kato76} (whose terminology
%and notations we follow here), the term ``positive'' means
%\((\mH x,x) > 0\) for any nonvanishing \(x \in \d (\mH)\).
This implies that all fractional powers of \(\mH\) and \(\mM\)
are also positive. Neither of these operators need be bounded (or have
bounded inverse).
\begin{theorem}\label{main}
Let \(\mH\) and \(\mM\) be positive selfadjoint operators in a Hilbert
space \(\x\) having the following property (A):
 \(\d (\mH^{1/2}) =  \d (\mM^{1/2})\) and the
norms \(\|\mH^{1/2}\cdot\|\) and  \(\|\mM^{1/2}\cdot\|\) are topologically
equivalent. Then the same property is shared by
\(\mH^{1/2}\) and  \(\mM^{1/2}\). The operators
\[
\overline{\mM^{-1/2}\mH^{1/2}},\quad  \overline{\mM^{1/2}\mH^{-1/2}},\quad
\overline{\mM^{-1/4}\mH^{1/4}},\quad  \overline{\mM^{1/4}\mH^{-1/4}},\quad
\]
\begin{equation}\label{allAB}
\overline{\mH^{-1/2}\mM^{1/2}},\quad  \overline{\mH^{1/2}\mM^{-1/2}},\quad
\overline{\mH^{-1/4}\mM^{1/4}},\quad  \overline{\mH^{1/4}\mM^{-1/4}}
\end{equation}
are well defined and bounded. Let
\begin{equation}\label{Tgeneral}
T = \overline{\mM^{1/2}\mH^{-1/2}} - \overline{\mM^{-1/2}\mH^{1/2}}
\end{equation}
and
\begin{equation}\label{Xgeneral}
X = \overline{\mM^{1/4}\mH^{-1/4}} - \overline{\mM^{-1/4}\mH^{1/4}}
\end{equation}
then
\begin{equation}\label{sqrtgeneral}
\|X\| \leq \frac{1}{2}\|T\|.
\end{equation}
\end{theorem}
\begin{proof}
The fact that the square roots inherit the property (A)
is  a consequence of L\"{o}wner type theorems
(see e.g.~ \cite{Kato76}, Ch.V, Th. 4.12). The corresponding pairs of operators
in (\ref{allAB}) are mutually adjoint e.g.
\(\overline{\mM^{-1/2}\mH^{1/2}}^* = \overline{\mH^{1/2}\mM^{-1/2}}\)
etc. Obviously, (\ref{Tgeneral}) and (\ref{Xgeneral})
reduce to \(T,X\) from (\ref{Tgeneral}), if the space is finite
dimensional. The equation (\ref{oursylvmatr}) becomes here
\begin{equation}\label{oursylvgeneral}
(X\mH^{-1/4}u,\mM^{1/4}v) + (X\mH^{1/4}u,\mM^{-1/4}v) = (Tu,v)
\end{equation}
for \(u \in \d _A = \d (\mH^{1/4})\cap\d (\mH^{-1/4})\)
and similarly for \(v\) and \(\mM\). We will now prove this.

The  left hand side of (\ref{oursylvgeneral}) equals
\[
(\overline{\mM^{1/4}\mH^{-1/4}}\mH^{-1/4}u,\mM^{1/4}v) -
(\mH^{-1/4}u,\overline{\mH^{1/4}\mM^{-1/4}}\mM^{1/4}v)
\]
\[
+ (\mH^{1/4}u,\overline{\mH^{-1/4}\mM^{1/4}}\mM^{-1/4}v) -
(\overline{\mM^{-1/4}\mH^{1/4}}\mH^{1/4}u,\mM^{-1/4}v) =
\]
\[
(\mH^{-1/2}u,\mM^{1/2}v) -
(u,v)
+ (u,v) -
(\mH^{1/2}u,\mM^{-1/2}v) =
\]
\[
(\mM^{1/2}\mH^{-1/2}u,v) - (\mM^{-1/2}\mH^{1/2}u,v) = (Tu,v).
\]
Now, substitute in (\ref{oursylvgeneral})
\begin{equation}\label{uvphipsi}
v = e^{-\mM^{-1/2}t}\mM^{-1/4}\phi,\quad
u = e^{-\mH^{-1/2}t}\mH^{-1/4}\psi
\end{equation}
for any \(\phi \in \d (\mM^{-1/2}),\ \psi \in \d (\mH^{-1/2})\).
Note that $$
\mM^{-1/4}\d (\mM^{-1/2})\qquad \text{and}\qquad
\mH^{-1/4}\d(\mH^{-1/2})$$
 are invariant under \(e^{-\mM^{-1/2}t}\),
\(e^{-\mH^{-1/2}t}\), respectively so, in (\ref{uvphipsi}) we have
\(u \in \d _A\) and \(v \in \d _\mM\). Then integrate
(\ref{uvphipsi}) and use partial integration:
\[
\int_0^s(Xe^{-\mH^{-1/2}t}\mH^{-1/2}\psi,e^{-\mM^{-1/2}t}\phi)dt +
\int_0^s(Xe^{-\mH^{-1/2}t}\psi,e^{-\mM^{-1/2}t}\mM^{-1/2}\phi)dt =
\]
\[
-\int_0^s(X\frac{d}{dt}e^{-\mH^{-1/2}t}\psi,e^{-\mM^{-1/2}t}\phi)dt +
\int_0^s(Xe^{-\mH^{-1/2}t}\psi,e^{-\mM^{-1/2}t}\mM^{-1/2}\phi)dt =
\]
\[
(X\psi,\phi) - (Xe^{-\mH^{-1/2}s}\psi,e^{-\mM^{-1/2}s}\phi)
+ \int_0^s(Xe^{-\mH^{-1/2}t}\psi,(-e^{-\mM^{-1/2}t}\mM^{-1/2})\phi)dt
\]
\[
+\int_0^s(Xe^{-\mH^{-1/2}t}\psi,e^{-\mM^{-1/2}t}\mM^{-1/2}\phi)dt =
\]
\[
(X\psi,\phi) - (Xe^{-\mH^{-1/2}s}\psi,e^{-\mM^{-1/2}s}\phi) =
\]
\[
 \int_0^s(Te^{-\mH^{-1/2}t}\mH^{-1/4}\psi,e^{-\mM^{-1/2}t}\mM^{-1/4}\phi)dt.
\]
In the limit \(s \to \infty\) by using the functional calculus for
\(\mH\), \(\mM\), respectively and monotone convergence for spectral integrals
we obtain
\[
e^{-\mH^{-1/2}s}\psi \to 0,\quad e^{-\mM^{-1/2}s}\phi \to 0
\]
in the norm. Hence
\begin{equation}\label{Xpsiphi}
(X\psi,\phi)
= \int_0^\infty(Te^{-\mH^{-1/2}t}\mH^{-1/4}\psi,e^{-\mM^{-1/2}t}\mM^{-1/4}\phi)dt
\end{equation}
where the integral on the right hand side is, in fact, Lebesgue
as shows the chain of inequalities in (\ref{estimation}) which are
valid in this general case as well. Here the identity
(\ref{Cidentity}) is used in the weak sense:
\[
\int_0^\infty(e^{-2Ct}C\phi,\phi)dt = (\phi,\phi)/2, \quad
\phi \in \d (C)
\]
for any positive selfadjoint \(C\). Thus,
\[
|(X\psi,\phi)|^2 \leq \|T\|^2(\psi,\psi)(\phi,\phi)/4.
\]
\end{proof}

\begin{remark}\label{forms}
The main assertion (\ref{sqrtgeneral}) of Theorem \ref{main} is
obviously equivalent to the following statement:
the inequality
\[
|m(\phi,\psi) - h(\phi,\psi)| \leq \varepsilon \sqrt{h(\phi,\phi)m(\phi,\psi)}
\]
implies
\[
|m_2(\phi,\psi) - h_2(\phi,\psi)|
\leq \frac{\varepsilon}{2} \sqrt{h_2(\phi,\phi)m_2(\psi,\psi)}
\]
where the sesquilinear forms \(h, m, h_2, m_2\) belong to the operators
$\mH$, $\mM$, $\mH^{1/2}$, $\mM^{1/2}$, respectively. Thus, our theorem will be
directly applicable to differential operators given in weak form.
\end{remark}

\begin{example}
Let again $\mH$ and $\mM$ be as in Example \ref{example:measur}. That is to say
take \(\mH\), \(\mM\) as selfadjoint realizations of the differential
operators
\[
-\frac{\partial}{\partial x}\alpha(x)\frac{\partial}{\partial x},\quad
-\frac{\partial}{\partial x}\beta(x)\frac{\partial}{\partial x},
\]
in the Hilbert space \(\H=L^2(I)\) (again  \(I\) can be a finite or infinite
interval) with the Dirichlet  boundary conditions and non-negative
bounded measurable functions $\alpha(x)$, $\beta(x)$ which satisfy
\[
|\beta(x) - \alpha(x)| \leq \varepsilon\sqrt{\beta(x)\alpha(x)}
\]
Now
\[
|(\mM^{1/2}\phi,\mM^{1/2}\psi) - (\mH^{1/2}\phi,\mH^{1/2}\psi)|^2 \leq
\left(\int_I |\beta(x) - \alpha(x)||\psi'(x)\phi'(x)|dx\right)^2 \leq
\]
\[
\varepsilon^2\int_I \alpha(x)|\psi'(x)|^2dt\int_I \beta(x)|\phi'(x)|^2dt
= \varepsilon^2\|\mH^{1/2}\phi\|^2\|\mM^{1/2}\psi\|^2
\]
hence \(\|T\| \leq \varepsilon\) and Theorem \ref{main}
applies yielding
\[
|(\mM^{1/4}\phi,\mM^{1/4}\psi) - (\mH^{1/4}\phi,\mH^{1/4}\psi)|  \leq
\frac{\varepsilon}{2}\|\mH^{1/4}\phi\|\|\mM^{1/4}\psi\|
\]
or, equivalently, in the terms as in Remark \ref{forms}
\[
|m_2(\phi,\psi) - h_2(\phi,\psi)|
\leq \frac{\varepsilon}{2} \sqrt{h_2(\phi,\phi)m_2(\psi,\psi)}.
\]
\end{example}

\section{Conclusion}
With this work we complete our study of the weak Sylvester equation which started
in \cite{GruVes02}. A notion of a weak Sylvester equation was introduced in
\cite{GruVes02} as a tool on a way to obtain invariant subspace estimates
for unbounded positive definite operators. With this paper we show that there
are applications of the concept of a weak Sylvester equation outside the theory
of Rayleigh--Ritz spectral approximations. We have extended out theory to
infinite dimensional invariant subspaces and have obtained estimates of
the difference between the corresponding spectral projections in all unitary invariant norms.
With this results we have developed a counterpart of the $\sin\Theta$
theorems from \cite{BathDavMcInt85} for perturbations of
operators which are only defined as quadratic forms.

Due to the very singular nature of integral representations (which can not be avoided by reformulation of the integrals)
of the solution to the equation (\ref{prvo2:e_slabaS1}), cf. formula (\ref{prvo2:e_integralnaR}), we were
not able to extend the technique from \cite{BathDavMcInt85} to
prove that in the setting of Theorem \ref{cor:unitary_invariant_2} assumption
$\tripleb F\tripleb<\infty$ also implies that there exists a bounded solution $T$ such
that $\tripleb T\tripleb<\infty$. We believe that this statement is true, but the proof
will have to remain a task for the future and will most likely require another technique.

An application of the concept to a perturbation of the square root of a
positive definite operator shows that there are other application areas for
weakly formulated operator equations and that the developed techniques are
(and hopefully will be) easily adaptable to new situations. The applications
which we have reported in this paper are presented as an illustration only.
Further applications will be the subject of the future work, cf. \cite{GruPhd}.
\bibliographystyle{abbrv}
%\bibliography{../bibliografija}
\def\cprime{$'$} \def\cprime{$'$} \def\cprime{$'$}

\end{document}